\documentclass[11pt]{amsart}

\usepackage{graphicx,epsfig}
\usepackage[dvipsnames]{xcolor}
\usepackage{amsfonts,amsmath,amssymb,amsthm}
\usepackage[normalem]{ulem}
\usepackage{comment,constants}
\usepackage{hyperref}
\usepackage{mathtools}
\usepackage{multirow}

\newtheorem{thm}{Theorem}[section]
\newtheorem{prop}[thm]{Proposition}
\newtheorem{lem}[thm]{Lemma}
\newtheorem{cor}[thm]{Corollary}
\newtheorem{conj}[thm]{Conjecture}
\newtheorem*{asm*}{Assumptions}
\newtheorem{asm}{Assumption}

\theoremstyle{remark}
\newtheorem{rem}[thm]{Remark}
\newtheorem*{rem*}{Remark}

\theoremstyle{definition}

\newcommand{\ra}{\rightarrow}

\newcommand{\N}{\mathbb N}     
\newcommand{\R}{\mathbb R}     
\newcommand{\Z}{\mathbb Z}     

\renewcommand{\a}{\alpha}
\renewcommand{\b}{\beta}

\renewcommand{\d}{\delta}

\newcommand{\e}{\varepsilon}

\renewcommand{\l}{\lambda}

\newcommand{\s}{\sigma}
\renewcommand{\th}{\theta}
\renewcommand{\epsilon}{\varepsilon}

\newcommand{\fl}[1]{\lfloor #1 \rfloor}  
\newcommand{\ind}[1]{ \mathbf{1}_{ \{ #1 \} } } 

\newconstantfamily{c}{symbol=c}
\newconstantfamily{CU}{symbol=C_U}

\newcommand{\w}{\omega}              
\renewcommand{\P}{\mathbb{P}}        
\newcommand{\E}{\mathbb{E}}          
 
\DeclareMathOperator{\Var}{Var}

\title[Recurrent excited random walks]{Functional limit laws for recurrent excited random walks with periodic cookie stacks}

\author{Elena Kosygina}
\address{Elena Kosygina\\One Bernard Baruch Way \\ Department of Mathematics, Box B6-230 \\ Baruch College \\ New York, NY 10010 \\ USA}
\email{elena.kosygina@baruch.cuny.edu}
\urladdr{http://www.baruch.cuny.edu/math/elenak/}
\thanks{E.\,Kosygina is partially supported by the Simons Foundation through a Collaboration Grant for Mathematicians \#209493.}

\author{Jonathon Peterson}
\address{Jonathon Peterson\\Purdue University\\Department of Mathematics\\150 N University Street\\West Lafayette, IN  47907\\USA}
\email{peterson@purdue.edu}
\urladdr{http://www.math.purdue.edu/~peterson}
\thanks{J. Peterson was partially supported by NSA grant H98230-15-1-0049.}


\begin{document}

\begin{abstract}
  We consider one-dimensional excited random walks (ERWs) with
  periodic cookie stacks in the recurrent regime.  We prove functional
  limit theorems for these walks which extend the previous results 
  in \cite{dkSLRERW}
  for excited random walks with ``boundedly
  many cookies per site.''  In particular, in the non-boundary
  recurrent case the rescaled excited random walk converges in the
  standard Skorokhod topology to a Brownian motion perturbed at
  its extrema (BMPE).  While BMPE is a natural limiting
  object for excited random walks with boundedly many cookies per
  site, it is far from obvious why the same should be true for our
  model which allows for infinitely many ``cookies'' at each site.
  Moreover, a BMPE has two parameters $\alpha,\beta<1$ and the scaling
  limits in this paper cover a larger variety of choices for $\a$ and
  $\b$ than can be obtained for ERWs with boundedly many cookies per
  site.

\end{abstract}

\maketitle

\section{Introduction}

Excited random walks (ERWs), also sometimes called cookie random walks, are self-interacting random walks where the transition probabilities of the walk depend on the local time of the walk at the current site. More precisely, a \emph{cookie environment} $\w = \{\w_x(j)\}_{x\in \Z, \, j\geq 1}$ is an element of $[0,1]^{\Z\times \N}$. For any fixed cookie environment $\w \in \Omega$, an ERW is a nearest-neighbor path $\{X_n\}_{n\geq 0}$ starting at $X_0 = 0$ and evolving so that on the $j$-th visit to a site $y \in \Z$ the walk moves right (resp.\ left) on the next step with probability $\w_y(j)$ (resp.\ $1-\w_y(j)$). 
That is, $P_\w$ is the law on nearest-neighbor paths on $\Z$ with $P_\w(X_0 = 0)$ and 
\begin{align*}
 P_\w\left( X_{n+1} = X_n+1 \, | \, X_0,X_1,\dots,X_n \right)
&= 1 - P_\w\left( X_{n+1} = X_n - 1 \, | \, X_0,X_1,\dots,X_n \right) \\
&= \w_{X_n}\left(\sum_{k=0}^n \ind{X_k=X_n} \right). 
\end{align*}
\begin{rem}
 The \emph{cookie} terminology comes from the following interpretation. One imagines a stack of cookies at each site in $\Z$. On each visit to a site the walker eats the next cookie in the stack at that site and the cookie creates an ``excitement'' that determines the transition probability of the next step. 
\end{rem}

The description of ERW given above was
for a fixed cookie environment $\w$, but one can also allow the cookie
environment to be random.  The probability distributions $P_\w$
defined above are called the \emph{quenched} laws of the ERW, while if
$\P$ is a probability distribution on the space of cookie environments
$\Omega$ then the \emph{averaged} law of the ERW is given by
\[
 P( \cdot ) = \int_\Omega P_\w(\cdot) \, \P(d\w). 
\]
To obtain some spatial regularity it is usually assumed that
under $\P$ the cookie stacks are stationary and ergodic under
  the shifts on $\Z$ (see, for example, \cite{zMERW}) or i.i.d. (the
  most common
  assumption). 
 That is, if $\w_x = \{\w_x(j)\}_{j\geq 1}$
 denotes the cookies stack at $x\in \Z$, then it is assumed that the
 sequence $\{\w_x\}_{x\in \Z}$ is either ergodic or i.i.d..

The recent review article \cite{kzERWsurvey} gives an
  extensive summary of many of the known results for ERW on $\Z^d$, $d\ge 1$. Here we will give a shorter summary of the
    known results for one-dimensional ERW that are relevant for the
    present paper.  The most extensive results for one-dimensional ERW
    are under the following assumptions on the cookie environments.
\begin{itemize}
 \item The cookie stacks are (spatially) i.i.d.
 \item There is an $M<\infty$ such that $\w_x(j) = 1/2$ for all $j>M$. 
\end{itemize}
We will refer to ERW under these assumptions as the
case of \emph{boundedly many cookies per site} since only the first $M$
cookies at each site give a non-zero drift (an ``excitement'' to the
right or left).  
For this model of ERW
many results are known including (but not limited
to) explicit criterion for recurrence/transience, a law of large
numbers with an explicit criteria for ballisticity, limiting
distributions, large deviation asymptotics, and scaling limits of the
occupation times of the right and left semi-axes
\cite{zMERW,bsCRWspeed,bsRGCRW,kzPNERW,kmLLCRW,dkSLRERW,pERWLDP,kzEERW,pXSDERW}.
In these results, many aspects of the behavior of the walk are
determined by a single explicit parameter,
\begin{equation}\label{deltadef}
 \d := \E\left[ \sum_{j\geq 1} (2\w_0(j)-1) \right],
\end{equation}
which is the expected total drift contained in the cookie stack at a fixed site. 
For instance, the walk is recurrent if and only if $\d \in [-1,1]$ \cite{zMERW,kzPNERW}, and the type of the limiting distributions are determined by the value of $\d$ \cite{bsRGCRW,kzPNERW,kmLLCRW,dkSLRERW}.

In this paper, instead of assuming boundedly many cookies per site, we
will consider the model of ERW with \emph{periodic cookie stacks}
which was first introduced in \cite{kosERWPC}. 
\begin{asm}\label{asm:per}
For given $N\in\N$ and $p_1,p_2,\ldots,p_N \in (0,1)$ with $\bar{p} = \frac{1}{N} \sum_{j=1}^N p_j = \frac{1}{2}$ every environment  $\w = \{\w_x(j)\}_{x\in \Z,j\ge 1}$ satisfies
\[
 \w_x(kN+j) = p_j, \quad \forall x \in \Z,\ k\geq 0, \text{ and } j=1,2,\ldots,N. 
\] 
\end{asm}
One can also consider ERW with periodic cookie stacks as in Assumption
\ref{asm:per} but with $\bar{p} \ne 1/2$.  In this
case, it was shown in \cite{kosERWPC}  and \cite{kpERWMC} that the ERW
is transient with non-zero speed and with a Gaussian limiting
distribution under diffusive scaling. Assumption \ref{asm:per}
restricts us to the critical case $\bar{p} = 1/2$ where the asymptotic
behavior of the walk is much more delicate.  Since under Assumption
\ref{asm:per} the sums $\sum_{j\geq 1}^n (2\w_0(j) - 1)$ oscillate as
$n\ra\infty$, there is no obvious way to generalize the known results
for boundedly many cookies per site  which are expressed
  in terms of a single parameter $\d$ in \eqref{deltadef}.  However,
the following result from \cite{kosERWPC} gives an explicit criterion
for recurrence/transience of ERW.
\begin{thm}[\cite{kosERWPC}]\label{thm:rt}
 Let Assumption \ref{asm:per} hold, and let $\th$ and $\tilde\th$ be defined by 
\begin{equation}\label{thetaper}
 \theta = \frac{\sum_{j=1}^N \sum_{i=1}^j (1-p_j)(2p_i-1)}{2 \sum_{j=1}^N p_j(1-p_j)}
\quad\text{and}\quad
 \tilde\theta = \frac{\sum_{j=1}^N \sum_{i=1}^j p_j(1-2p_i)}{2 \sum_{j=1}^N p_j(1-p_j)}.
\end{equation}
\begin{enumerate}
 \item If $\th>1$ then $P(\lim_{n\ra\infty} X_n = \infty) = 1$. 
 \item If $\tilde\th>1$ then $P(\lim_{n\ra\infty} X_n = -\infty) = 1$. 
 \item If $\max\{\th,\tilde\th\} \leq 1$ then $P( \liminf_{n\ra\infty} X_n = -\infty, \, \limsup_{n\ra\infty} X_n = \infty) = 1$. 
\end{enumerate}
\end{thm}

\begin{rem}\label{rem:thsum}
The fact that the parameters $\th$ and $\tilde\th$ cannot both be greater than one follows from the relation $\th + \tilde\th = 1 - \frac{N}{4\sum_{j=1}^N p_j(1-p_j)}.$ This identity can be obtained either from the formulas in \eqref{thetaper} and some algebra or as a consequence of a more general argument in \cite[Proposition 4.3]{kpERWMC}. 
\end{rem}

\begin{rem}
  The proof of Theorem \ref{thm:rt} in \cite{kosERWPC} uses an
  approach based on Lyapounov functions to give criteria for
  recurrence and transience.  Another proof of Theorem \ref{thm:rt}
  was given in \cite{kpERWMC} for a more general model of ERW where
  the cookie stacks at each site come from independent realizations of
  a finite state Markov chain.  This more general model includes both
  periodic cookie stacks and certain models of bounded cookie stacks
  as special cases.  The proof in \cite{kpERWMC} relied on certain
  tail asymptotics for regeneration times of a related Markov chain.
  This method had the advantage of also leading to further results
  such as a criterion for ballisticity and a characterization of the
  limiting behavior in the transient cases; results which were
  previously only known for the case of bounded cookie stacks.  Not
  covered in \cite{kpERWMC} were the scaling limits of ERWs in the
  recurrent cases. This is the topic of the current
  paper.
\end{rem}

\begin{rem}
While for ERW with bounded cookie stacks the recurrence/transience, ballisticity, and limiting distributions in the transient cases depend only on the single parameter $\d$ defined in \eqref{deltadef}, the more general results in \cite{kosERWPC} and \cite{kpERWMC} for ERW with periodic (or Markovian) cookie stacks rely on two parameters $\theta$ and $\tilde\th$. 
In the special case of bounded cookie stacks these parameters are $\theta = \d$ and $\tilde\theta = -\d$, but in general it is not the case that $\theta + \tilde\theta = 0$ (see Remark \ref{rem:thsum} above). 
\end{rem}

\subsection{Main results: functional limit theorems in the recurrent regime}

To review the known results for recurrent ERW with boundedly many cookies per site, we first must recall the definition of a \emph{perturbed Brownian motion.} For fixed parameters $\a,\b \in (-\infty,1)$, a $(\a,\b)$-perturbed Brownian motion is a solution 
$Z^{\a,\b}_\cdot$ 
to the functional equation 
\begin{equation}\label{pBM}
Z^{\a,\b}_0 = 0 \quad\text{and}\quad 
Z^{\a,\b}_t = B_t + \a \sup_{s\leq t} Z^{\a,\b}_s + \b \inf_{s \leq t} Z^{\a,\b}_s, 
\quad \text{for } t > 0,
\end{equation}
where $B_t$ is a standard Brownian motion. 
It was shown in \cite{pwPBM,cdPUPBM} that if $\a,\b<1$ then there is almost surely a pathwise unique solution of \eqref{pBM} that is continuous and adapted to the filtration of the Brownian motion.\footnote{A perturbed Brownian motion does not exist if $\alpha\ge 1$ or $\beta\ge 1$.}
The following functional limit theorems were proved in \cite{dkSLRERW} for recurrent ERW with boundedly many cookies per site. 
\begin{itemize}
 \item \textbf{Boundary case.} If $\d = 1$ then there exists a constant $a>0$ such that $\{\frac{X_{\fl{nt}}}{a \sqrt{n}( \log n)} \}_{t\geq 0}$ converges in distribution to the running maximum of a Brownian motion $B^*(t) = \sup_{s\leq t} B_s$. Similarly, if $\d=-1$ then the rescaled ERW converges to the running minimum of a Brownian motion. 
 \item \textbf{Non-boundary case.} If $\d \in (-1,1)$, then $\{ \frac{X_{\fl{nt}}}{\sqrt{n}} \}_{t\geq 0}$ converges in distribution to a $(\d,-\d)$-perturbed Brownian motion. 
\end{itemize}
Our main results are similar functional limit theorems for ERW with periodic cookie stacks. 
Here, and throughout the paper, $D([0,\infty))$ will denote the space of c\'adl\'ag functions equipped with the Skorokhod $J_1$ topology, and convergence in distribution on this space will be denoted by $\overset{J_1}{\Longrightarrow}$.

\begin{thm}[Recurrent ERW - boundary case]\label{thm:boundary}
Let Assumption \ref{asm:per} hold.
If $\th = 1$ then there exists a constant $a>0$ such that 
\[
 \left\{  \frac{X_{\fl{n t}}}{a \sqrt{n}(\log n)} \right\}_{t\geq 0} \xRightarrow[n\ra\infty]{J_1} \left\{ B^*_t \right\}_{t\geq 0}, 
\]
where $B^*_t = \sup_{s\leq t} B_s$ is the running maximum of a standard Brownian motion. 
If $\tilde\th = 1$ a similar scaling limit holds with the limiting process instead being the running minimum of a Brownian motion. 
\end{thm}

\begin{thm}[Recurrent ERW - non-boundary case]\label{thm:pbm}
 Let Assumption \ref{asm:per} hold. If $\max\{\theta,\tilde\theta\} < 1$, then 
\[
 \left\{ \frac{X_{\fl{nt}}}{a \sqrt{n}} \right\}_{t\geq 0} \xRightarrow[n\ra\infty]{J_1} \left\{ Z^{\theta,\tilde\theta}_t \right\}_{t\geq 0}, \qquad \text{where } a = \frac{1}{2} \left( \frac{1}{N} \sum_{i=1}^N p_i(1-p_i) \right)^{-1/2}, 
\]
and $Z^{\theta,\tilde\theta}$ is a $(\theta,\tilde\theta)$-perturbed Brownian motion as defined in \eqref{pBM}. 
\end{thm}

The proof of Theorem \ref{thm:boundary} is follows word
  for word the proof of of \cite[Theorem 1.2]{dkSLRERW} for boundedly
  many cookies per site once we substitute the necessary tail
  asymptotic results for the associated branching-like processes
  (proved in \cite{kpERWMC} and re-stated in Theorem \ref{thm:BLPtail}
  below) for the corresponding results in the case of boundedly many
  cookies per site.  We will therefore omit the proof of Theorem
  \ref{thm:boundary} and focus on the proof of Theorem \ref{thm:pbm}.

For recurrent ERW with boundedly many cookies per site it is easy to
see why the scaling limit would be a perturbed Brownian motion. After
a large number of steps, one can expect that the walk will have
visited the sites in the interior of its range a large number of
times. If the walk only experiences ``excitement'' in the first $M$
visits to a site then it is intuitively obvious that the limiting
process should behave like a Brownian motion when it is away from its
running minimum or maximum and should experience some additional drift
at the edge of its current range.
For ERW with periodic cookie stacks it is not nearly so obvious why
the scaling limit should be a Brownian motion in the interior of its
range.  In fact, while the ERW does scale to a Brownian motion in the
interior of the range, since the scaling parameter $a$ in Theorem
\ref{thm:pbm} is larger than one\footnote{except in the simple random walk case with $p_i=1/2$ for all $i$.} it is evident that consecutive steps of the ERW
in the interior of the range are  quite strongly correlated.

\subsection{Overview of ideas and structure of the paper}
Many of the recent results about one-dimensional ERWs
  rely on the analysis of certain \emph{branching-like processes}
which are related to the directed edge local times of the random
walk. Similar ideas have been used previously in the study of other
non-standard random walks
\cite{kksStable,tTSAWGBR,tTSAW,tGRK,tLTWRRW}.  For ERW,
these methods were first used in \cite{bsCRWspeed} and have since
become the primary tool for the study of one-dimensional ERWs.

The connection between random walks and branching-like processes is in
the spirit of the Ray-Knight theorems which relate the local times of
a one-dimensional Brownian motion with certain squared Bessel
processes.  B.\,T\'oth used this connection to prove
  generalized Ray-Knight theorems for a large class of
  self-interacting random walks
  \cite{tTSAWGBR,tTSAW,tGRK,tLTWRRW}. These theorems were a key tool in
  obtaining limiting distributions for the random walk stopped at an
  independent exponential random time. Limiting distributions at
  deterministic late times still remain an open problem. We note, that
  for a certain sub-class of the self-interacting random walks
  B.\,T\'oth identified the limiting distributions with
  the one-dimensional marginal distributions of Brownian motion
  perturbed at its extrema \cite[Remark on p. 1334]{tGRK}.

  More recently, similar Ray-Knight theorems for ERWs have been proven
  and used in \cite{kmLLCRW,kzEERW,DK14,kpERWMC}.  In the present
  paper, we are able to combine some of the ideas introduced by T\'oth
  together with a martingale decomposition of the excited random walk
  as in \cite{dCLTERW, dkSLRERW} to prove a full
  process-level convergence to Brownian motion perturbed at its
  extrema.

The remainder of the paper is structured as follows. 
In Section \ref{sec:BLP} we recall the definitions of the branching-like processes associated with the ERWs.   
We will review the construction of these branching-like processes, their connection with the directed-edge local times of the ERW, and some results from \cite{kpERWMC} regarding tail asymptotics and scaling limits of these processes. 
In Section \ref{sec:prelim} we will prove some preliminary results in
preparation for the proof of Theorem \ref{thm:pbm}. In particular,
using the connection with the branching-like processes we will show
that diffusive scaling is the right scaling to obtain a limiting
distribution and that for any fixed $\gamma\in(0,1/2)$ and a
sufficiently large time $n$ most of the sites in the interior of the
range have been visited at least $n^{\gamma}$ times.  Finally, in
Section \ref{sec:pbm} we give the proof of Theorem \ref{thm:pbm}.  The
key to the proof is Lemma \ref{lem:Cnapprox} which gives sufficient
control on the total drift contained in the ``cookies'' used by the
ERW in the first $n$ steps. Again the connection with
the branching-like processes is essential.  We close the paper in
Section \ref{sec:Markov} with a brief discussion of the more general
model of ERW with Markovian cookie stacks introduced in \cite{kpERWMC}
and explain the difficulty in extending Theorem \ref{thm:pbm} to this
more general model.

\section{Related branching-like processes}\label{sec:BLP}

In this section we will introduce four Markov chains which we will
refer to as ``branching-like processes'' (BLPs) that are related to
the directed-edge local times of the random walk.  
We will also recall some important results concerning the BLPs proved in \cite{kpERWMC}; in particular, we will recall certain tail asymptotic results (Theorem \ref{thm:BLPtail}) and  the fact that scaling limits of BLPs are squared Bessel processes (Theorem \ref{thm:da}).

\subsection{Construction of the BLP}

To prepare for both the construction of the BLPs and the connection with the random walk, we first recall the following simple construction of the ERW. 
While we are primarily interested in ERW with periodic cookie stacks in this paper, the results of this section are more general (with the exception of the specific formulas for parameters in Section \ref{sssec:theta} below), and thus we will give the construction of the BLPs in the more general setting of random cookie environments that are (spatially) i.i.d.\, (that is $\{\w_x\}_{x\in\Z}$ is i.i.d.\ under the distribution $\P$ on cookie environments).

Given an environment $\w = \{\w_x(j)\}_{x\in\Z,\, j\geq 1}$ we let $\{\xi_x(j)\}$ be a family of independent Bernoulli random variables with $\xi_x(j) \sim \text{Ber}(\w_x(j))$. Then the path $\{X_n\}_{n\geq 0}$ of the ERW can be constructed iteratively as follows. If $X_n = x$ and $\sum_{k=0}^n \ind{X_k = x} = j$, then $X_{n+1} = x + (2\xi_x(j) - 1)$. 
That is, upon visiting a site $x$ for the $j$-th time the walk steps to the right if $\xi_x(j) = 1$ and to the left if $\xi_x(j) = 0$. 

We will now use these Bernoulli random variables $\xi_x(j)$ to construct the BLPs. For this construction we will only need to consider the sequence $\{\xi_x(j)\}_{j\geq 1}$ for a fixed $x$, and since these have the same distribution for each $x$ (under the averaged measure) we will, for simplicity of notation, simply use $\{\xi(j)\}_{j\geq 1}$ to denote one such sequence. 
Next, for any $m\geq 0$ let 
\[
 S_m = \inf\left\{ k\geq 0: \sum_{j=1}^{k+m}(1-\xi(j)) = m \right\} \quad\text{and}\quad F_m = \inf\left\{ k\geq 0: \sum_{j=1}^{k+m} \xi(j) = m \right\}. 
\]
Note that by the convention that an empty sum is equal to zero, we have that $S_0 = 0$ and $F_0 = 0$. 
If we refer to a Bernoulli random variable $\xi(j)$ as a ``success'' if $\xi(j) = 1$ and a ``failure'' if $\xi(j)=0$ then $S_m$ is the number of successes before the $m$-th failure and $F_m$ is the number of failures before the $m$-th success in the sequence of Bernoulli trials $\{\xi(j)\}_{j\geq 1}$. 
Having introduced this notation, we may now define the BLPs $U,\hat{U},V$ and $\hat{V}$ to be Markov chains on $\Z_+ = \{0,1,2,\ldots\}$ with the following transition probabilities. 
\begin{align*}
& P(U_{k+1} = n \, | \, U_k = m) = P( S_m = n ) & & P(V_{k+1} = n \, | \, V_k = m) = P(F_m = n) \\
& P(\hat{U}_{k+1} = n \, | \, \hat{U}_k = m) = P( S_{m+1} = n )  &&  P(\hat{V}_{k+1} = n \, | \, \hat{V}_k = m) = P( F_{m+1} = n )
\end{align*}
Note since $S_0=F_0=0$ that the BLPs $U$ and $V$ are absorbing at the state $x=0$. On the other hand, the Markov chains $\hat{U}$ and $\hat{V}$ are irreducible Markov chains. 

To explain the terminology of ``branching-like processes'' note that if $\w_x(j) \equiv \a \in (0,1)$ for all $j\geq 1$ (in this case the ERW is just a simple random walk) then the processes $U$ and $V$ are branching processes with offspring distributions that are Geo($1-\a$) and Geo($\a$), respectively, and the processes $\hat{U}$ and $\hat{V}$ are branching processes with the same offspring distributions but with an extra immigrant prior to reproduction in each generation. 
For ERW with boundedly many cookies per site, the processes $U,\hat{U},V$ and $\hat{V}$ can be interpreted as branching processes with migration (see \cite{bsCRWspeed,kzPNERW}). 
For a discussion of the branching-like structure of the processes in the more general case of Markovian (or periodic) cookie stacks see \cite[Section 2]{kpERWMC}.

\subsection{Connection with ERW}\label{connection}
Before studying properties of the above BLPs more in depth, we will first recall the connection of these processes with the directed edge local times of the random walk. To this end, let $\l_{x,m}$ be the stopping times for the ERW defined for $x \in \Z$ and $m\geq 0$
\[
 \l_{x,0} = \inf\{ n\geq 0: \, X_n = x \} \quad \text{and} \quad \l_{x,m} = \inf\{ n>\l_{x,m-1}: X_n = x \}. 
\]
That is, $\l_{x,m}$ is the time of the $(m+1)$-st visit to $x\in \Z$. Note that these stopping times could in theory be infinite, but since we are concerned in this paper only with recurrent ERW the stopping times $\l_{x,m}$ are almost surely finite. 
For a fixed $x \in \Z$ and $m\geq 0$, define the directed edge local time processes  $\mathcal{E}^{(x,m)}_y$ and $\mathcal{D}^{(x,m)}_y$ as follows.
\[
 \mathcal{E}^{(x,m)}_y = \sum_{k=0}^{\l_{x,m}-1} \ind{X_k = y, \, X_{k+1} = y+1}
\quad \text{and}\quad
 \mathcal{D}^{(x,m)}_y = \sum_{k=0}^{\l_{x,m}-1} \ind{X_k = y, \, X_{k+1} = y-1}. 
\]
That is, $\mathcal{E}_y^{(x,m)}$ and $\mathcal{D}^{(x,m)}_y$ give the number of steps to the right and left, respectively, from the site $y$ by time $\l_{x,m}$. 

For fixed $x\in \Z$ and $m\geq 0$, it follows from the construction of the ERW in terms of the Bernoulli random variables $\xi_y(j)$ that 
\begin{equation}\label{EDinitial}
 \mathcal{E}^{(x,m)}_x = \sum_{j=1}^m \xi_x(j) \quad\text{and}\quad \mathcal{D}^{(x,m)}_x = \sum_{j=1}^m (1-\xi_x(j)).
\end{equation}
We claim that the processes $\{\mathcal{E}^{(x,m)}_{x+k}\}_{k\geq 0}$ and $\{\mathcal{D}^{(x,m)}_{x-k}\}_{k\geq 0}$ are Markov chains with transition probabilities which are the same as the BLPs defined above. 
\begin{figure}[ht]
 \includegraphics[width=3in]{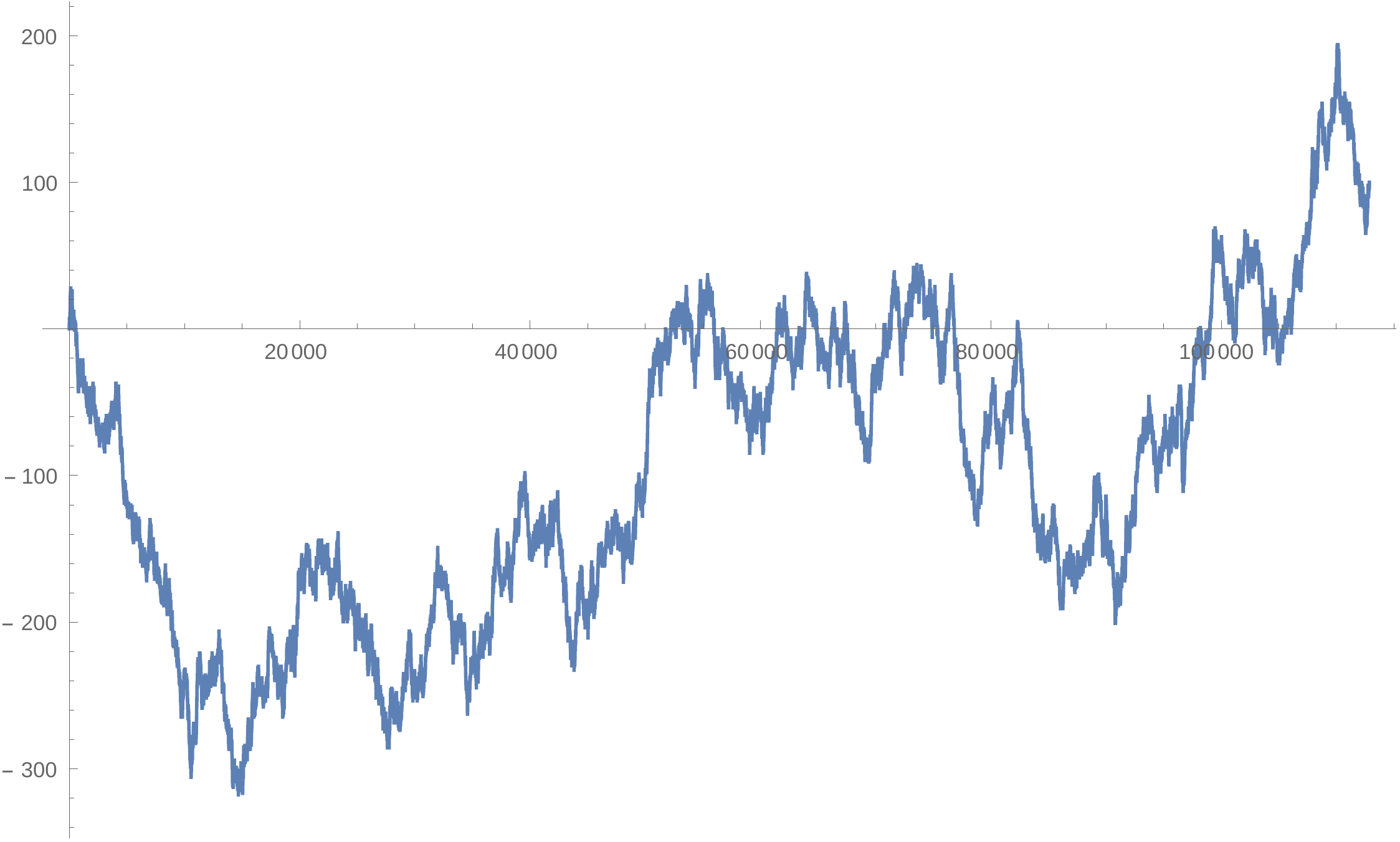} \hspace{.3in} \includegraphics[width=3in]{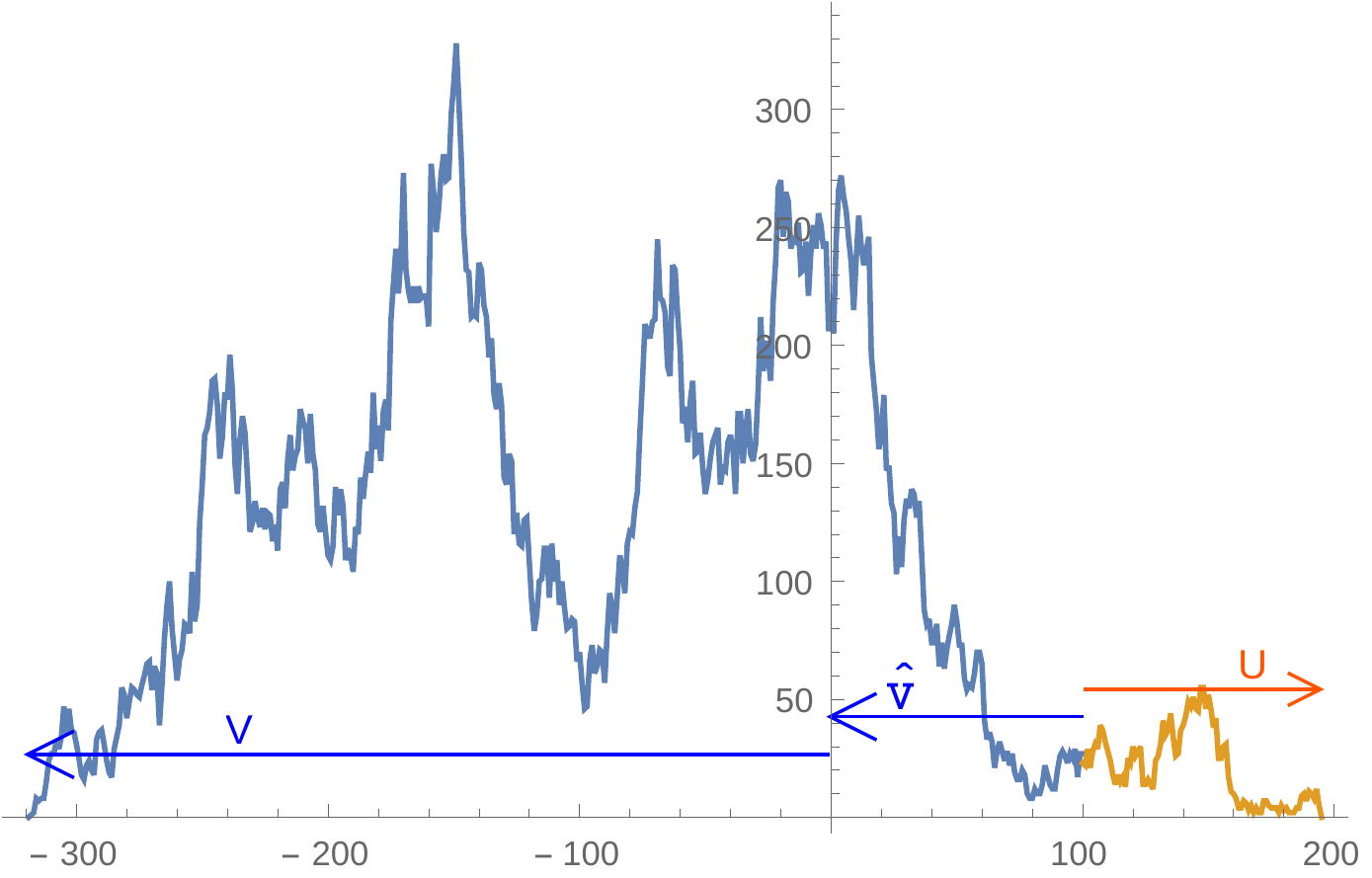}
\caption{On the left is a simulation of an ERW in a cookie with periodic cookie stacks of the form $\w_x = (0.7,0.3,0.7,0.3,\ldots)$ at each site $x\in\Z$ (in this case $\th = \frac{1}{7}$ and $\tilde\th = -\frac{1}{3}$). The path of the walk is simulated until the walk visits the site $x=100$ for the 51-st time. On the right are corresponding plots of the processes $\{\mathcal{D}_y^{(100,50)}\}_{y\leq 100}$ (in blue) and $\{\mathcal{E}_y^{(100,50)}\}_{y\geq 100}$ (in orange). The process in orange is a BLP of the form $U$. The process in blue (viewed from right to left) is a BLP which is of the form $\hat{V}$ to the right of the origin and $V$ to the left of the origin.} \label{fig:BLPplot}
\end{figure}
For specificity we consider first the case when $x>0$.
Since the random walk ends at site $x>0$ at time $\l_{x,m}$, it follows that if there are $\mathcal{E}^{(x,m)}_{x+k} = \ell$ jumps to the right from $x+k$ by time $\l_{x,m}$ then there are also $\ell$ jumps to the left from $x+k+1$ by time $\l_{x,m}$. Thus, $\mathcal{E}_{x+k+1}^{(x,m)}$ equals the number of jumps to the right from $x+k+1$ before the $\ell$-th jump to the left, or equivalently the number of successes before the $\ell$-th failure in the Bernoulli sequence $\{\xi_{x+k+1}(j)\}_{j\geq 1}$. Therefore, for any sequence $\ell_1,\ell_2,\cdots,\ell_{k+1} \in \Z_+$ we have 
\begin{align}
 P\left( \mathcal{E}_{x+k+1}^{(x,m)} = \ell_{k+1} \, \bigl| \,  \mathcal{E}_{x+i}^{(x,m)} = \ell_i, \, 1\leq i\leq k \right) 
&= P\left( \mathcal{E}_{x+k+1}^{(x,m)} = \ell_{k+1}  \, \bigl| \,  \mathcal{E}_{x+k}^{(x,m)} = \ell_k, \right) \nonumber \\
&= P\left( U_1 = \ell_{k+1} \, | \, U_0 = \ell_k \right).  \label{EMarkov}
\end{align}
That is, $\{\mathcal{E}_{x+k}^{(x,m)}\}_{k\geq 0}$ is a Markov chain with the same transition probabilities as the BLP $U$. 
(Note that for the first equality in \eqref{EMarkov} we use that the cookie environment $\w = \{\w_x\}_x$ is (spatially) i.i.d.)
The analysis of the process $\mathcal{D}^{(x,m)}_{x-k}$ when $x>0$ is similar, but slightly more complicated.
If there are $\ell$ steps to the left from $x-k$ by time $\l_{x,m}$ (i.e., $\mathcal{D}^{(x,m)}_{x-k} = \ell$) then the number of steps to the right from $x-k-1$ by time $\l_{x,m}$ is $\ell$ if $x-k-1<0$ and $\ell+1$ if $x-k-1 \geq 0$. 
This is because every jump from $x-k$ to $x-k-1$ is followed by a return from $x-k-1$ to $x-k$, but there is also an initial jump from $x-k-1$ to $x-k$ if $x-k-1\geq 0$. 
Therefore, similar to \eqref{EMarkov} we can conclude that 
\begin{align}
 P\left( \mathcal{D}_{x-k-1}^{(x,m)} = \ell_{k+1} \, \bigl| \,  \mathcal{D}_{x-i}^{(x,m)} = \ell_i, \, 1\leq i\leq k \right) 
&= P\left(  \mathcal{D}_{x-k-1}^{(x,m)} = \ell_{k+1} \, \bigl| \,  \mathcal{D}_{x-k}^{(x,m)} = \ell_k \right) \nonumber \\
&= \begin{cases}
    P\left( V_1 = \ell_{k+1} \, | \, V_0 = \ell_k \right) & x-k-1 < 0 \\ 
    P\left( \hat{V}_1 = \ell_{k+1} \, | \, \hat{V}_0 = \ell_k \right) & x-k-1 \geq 0.
   \end{cases}
\label{DMarkov}
\end{align}

The analysis of the directed edge local time processes is similar in the cases $x=0$ and $x<0$. A summary of the correspondence of the directed edge local time processes to the BLPs in the different cases is given in Table \ref{table:ED_BLP}. 
\begin{table}[ht]
\begin{tabular}{|l|l|c|}
\hline
 \textbf{Case} & \textbf{Directed edge local time} & \textbf{BLP} \\
\hline\hline
\multirow{3}{*}{$x<0$}
 & $\mathcal{E}^{(x,m)}_x, \mathcal{E}^{(x,m)}_{x+1},\ldots, \mathcal{E}^{(x,m)}_{-1},\mathcal{E}^{(x,m)}_0$ & $\hat{U}$ \\ \cline{2-3}
 & $\mathcal{E}^{(x,m)}_0, \mathcal{E}^{(x,m)}_{1},\mathcal{E}^{(x,m)}_{2},\ldots$ & $U$ \\ \cline{2-3}
 & $\mathcal{D}^{(x,m)}_x, \mathcal{D}^{(x,m)}_{x-1}, \mathcal{D}^{(x,m)}_{x-2},\ldots$ & $V$\\
\hline\hline
\multirow{2}{*}{$x=0$}
 & $\mathcal{E}^{(x,m)}_0, \mathcal{E}^{(x,m)}_{1},\mathcal{E}^{(x,m)}_{2},\ldots$ & $U$ \\ \cline{2-3}
 & $\mathcal{D}^{(x,m)}_0, \mathcal{D}^{(x,m)}_{-1}, \mathcal{D}^{(x,m)}_{-2},\ldots$ & $V$\\ 
\hline\hline
\multirow{3}{*}{$x>0$}
 & $\mathcal{E}^{(x,m)}_x, \mathcal{E}^{(x,m)}_{x+1},\mathcal{E}^{(x,m)}_{x+2},\ldots$ & $U$ \\ \cline{2-3}
 & $\mathcal{D}^{(x,m)}_x, \mathcal{D}^{(x,m)}_{x-1}, \ldots \mathcal{D}^{(x,m)}_1, \mathcal{D}^{(x,m)}_0$ & $\hat{V}$\\ \cline{2-3}
 & $\mathcal{D}^{(x,m)}_0, \mathcal{D}^{(x,m)}_{-1}, \mathcal{D}^{(x,m)}_{-2},\ldots$ & $V$\\
\hline
\end{tabular}
\caption{The directed edge local time processes are Markov chains with initial conditions given by \eqref{EDinitial} and transition probabilities corresponding to BLP as given in this table.}
\label{table:ED_BLP}
\end{table}

We note also that the directed edge local times can be used to represent the local time of the random walk at sites. That is, if $\mathcal{L}(n;x) = \sum_{k=0}^{n-1} \ind{X_k = x}$ is the number of visits to a site $x$ before time $n$, then 
\begin{equation}\label{LED}
 \mathcal{L}(\l_{x,m}; y) = \mathcal{D}_y^{(x,m)} + \mathcal{E}_y^{(x,m)}
= 
\begin{cases}
 m & y=x \\
\mathcal{D}_y^{(x,m)} + \mathcal{D}_{y+1}^{(x,m)} + \ind{0\leq y < x} & y < x \\
\mathcal{E}_{y-1}^{(x,m)} + \mathcal{E}_y^{(x,m)} + \ind{x<y \leq 0} & y > x. 
\end{cases}
\end{equation}
The first equality in \eqref{LED} is obvious since every visit to $y$ must result in a jump to the right or left. 
For the second equality in \eqref{LED} the formula in the case $y=x$ is clear by the definition of the stopping time $\l_{x,m}$. In the case $y>x$ the second equality in \eqref{LED} follows from the fact that $\mathcal{D}_y^{(x,m)} = \mathcal{E}_{y-1}^{(x,m)} + \ind{x<y\leq 0}$ since every jump to the left from $y$ can be paired with a preceding jump to the right from $y-1$ except in the case when $x<y\leq 0$ where there is no such corresponding jump for the first jump to the left from $y$.

\subsection{Previous results}
As noted above, the BLPs have been studied quite extensively in the case of ERW with boundedly many cookies per site, and many of these results were extended to the case of periodic (and even Markovian) cookie stacks in \cite{kpERWMC}. 
In this subsection we will recall some of these results which will be of importance in the current paper. 

\subsubsection{The parameters $\th$ and $\tilde\th$}\label{sssec:theta} We begin by recalling the connection of the parameters $\th$ and $\tilde\th$ defined in \eqref{thetaper} with the BLPs defined above. 
These parameters were defined in \cite{kosERWPC} in terms of asymptotics of the mean and variance of the BLPs when the BLPs are large. 
The parameters are given by $\th = \frac{2\rho}{\nu}$ and $\tilde\th = \frac{2\tilde\rho}{\nu}$ where 
\begin{equation}\label{BLPrhonu}
\begin{split}
 &\rho = \lim_{n\ra\infty} E[U_1 \, | \, U_0 = n] - n, \qquad \tilde\rho = \lim_{n\ra\infty} E[V_1 \, | \, V_0 = n] - n, \\
\text{and}\quad 
 &\nu = \lim_{n\ra\infty} \frac{\Var(U_1 \, | \, U_0 = n)}{n} = \lim_{n\ra\infty} \frac{\Var(V_1 \, | \, V_0 = n)}{n}. 
\end{split}
\end{equation}
In fact, it was shown in \cite{kosERWPC} that the convergence in the above definitions of $\rho$ and $\tilde\rho$ is exponentially fast in $n$. Since we will use this later, we record here that 
\begin{equation}\label{rhoroc}
 \left| E[U_1 \, | \, U_0 = n] - n - \rho \right| \leq C e^{-cn}. 
\end{equation}
for some constants $C,c>0$. 
The explicit formulas for $\th$ and $\tilde\th$ in \eqref{thetaper} follow from the following explicit formulas for $\rho$, $\tilde\rho$ and $\nu$ derived for periodic cookie stacks in \cite{kosERWPC}. 
\begin{equation}\label{rhonuforms}
 \rho = \frac{2}{N} \sum_{j=1}^N (1-p_j) \sum_{i=1}^j (2p_i-1), 
\quad 
\tilde\rho = \frac{2}{N} \sum_{j=1}^N p_j \sum_{i=1}^j (1-2p_i), 
\quad \text{and}\quad
\nu = \frac{8}{N} \sum_{i=1}^N p_i(1-p_i)\le 2. 
\end{equation}
Finally, we also note the following relation between the parameters $\rho$, $\tilde\rho$, and $\nu$ which will be used in the proof of Theorem \ref{thm:pbm} below. 
\begin{equation}\label{rhosums}
\rho + \tilde\rho = \frac{\nu}{2}-1. 
\end{equation}
This identity follows from the formulas in \eqref{rhonuforms} and the assumption that $\bar{p} = \frac{1}{N} \sum_{i=1}^N p_i = 1/2$. The analogs of \eqref{rhoroc}-\eqref{rhosums} for a more general model were obtained in \cite[Propositions 3.1, 4.3, and (37)-(38)]{kpERWMC}. 

\subsubsection{Tail asymptotics}
The relevance of the parameters $\th$ and $\tilde\th$ defined above is that they determine certain tail asymptotics for the BLPs. 
To provide some unified notation for studying these tail asymptotics we will adopt the following notation for certain hitting times of a stochastic process. 
If $\{Z_i\}_{i\geq 0}$ is a stochastic process, then for $x\in \R$ we will let 
\[
 \s_x^Z := \inf\{ i>0: Z_i \leq x \}. 
\]
\begin{thm}[Theorems 2.6 and 2.7 in \cite{kpERWMC}]\label{thm:BLPtail}
Let $Z$ be one of the BLPs $U, \hat{U}, V$, or $\hat{V}$, and let $s_Z$ be defined for each of these cases by 
\begin{equation}\label{gammaZ}
 s_U = 1-\th, \quad s_{\hat{U}} = \tilde\th, \quad s_V = 1-\tilde\th, \quad\text{and}\quad s_{\hat{V}} = \th. 
\end{equation}
Let $m\geq 1$ (or let $m\geq 0$ if $Z$ is either $\hat{U}$ or $\hat{V}$). 
\begin{enumerate}
 \item If $s_Z < 0$, then $P(\s_0^Z = \infty \, | \, Z_0 = m ) > 0$. 
 \item If $s_Z \geq 0$, then there exist constants $C^Z_1(m),C^Z_2(m)>0$ such that 
\[
 \lim_{n\ra\infty} n^{s_Z} P\left( \s_0^Z > n \, | \, Z_0 = m \right) = C^Z_1(m), 
\quad\text{and}\quad
 \lim_{n\ra\infty} n^{s_Z/2} P\left( \sum_{i=0}^{\s_0^Z-1} Z_i > n \, \biggl| \, Z_0 = m \right) = C^Z_2(m), 
\]
where for $s_Z=0$ we replace $n^{s_Z}$ with $\ln n$.
\end{enumerate}
\end{thm}

\subsubsection{Squared Bessel processes}
Finally, we recall the following  diffusive scaling limits for the BLPs. 
The asymptotics of the mean and variance of the BLPs in \eqref{BLPrhonu} suggest the following diffusion approximations for the BLPs which were proved in \cite{kpERWMC}. 

\begin{thm}[Lemma 6.1 in \cite{kpERWMC}]\label{thm:da}
Fix $y>\e>0$
and a sequence $z_n \in \Z_+$ with $z_n/n \ra y$ as $n\ra\infty$. Let $Z^n_\cdot$ be a sequence of one of the BLPs $U,\hat{U},V$, or $\hat{V}$ with initial conditions $Z^n_0 = z_n$, and let $Y^{\e,n}(t) = \frac{Z_{nt \wedge \s_{\e n}^Z}}{n}$ for $t\geq 0$. Then, with respect to the standard Skorokhod($J_1$) topology on the space of c\'adl\'ag functions on $[0,\infty)$ the process $Y^{\e,n}$ converges in distribution to $\{Y(t \cdot \wedge \s_\e^Y)\}_{t\geq 0}$, where $Y$ is the solution of 
\[
 dY(t) = b_Z \, dt + \sqrt{\nu \, Y(t)} \, dB(t), \quad Y(0) = y, 
\]
where $B(t)$ is a standard Brownian motion, and the drift constant $b_Z$ for the four possible BLPs is given by 
\[
 b_U = \rho, \quad
 b_{\hat{U}} = 1+\rho, \quad
 b_V = \tilde\rho, \quad\text{and}\quad
 b_{\hat{V}} = 1 + \tilde\rho. 
\]
\end{thm}
\begin{rem}
 Note that $2Y(t)$ is a (time-rescaled) squared Bessel process of generalized dimension $\frac{4 b_Z}{\nu}$. The exponents $s_Z$ defined in \eqref{gammaZ} which are used in Theorem \ref{thm:BLPtail} are related to this generalized dimension by the relation $s_Z = 1-\frac{1}{2}\left( \frac{4 b_Z}{\nu} \right)$. 
\end{rem}

\section{Preliminaries} \label{sec:prelim}


\subsection{Diffusive scaling}
We begin by proving some lemmas which indicate that the diffusive scaling is the correct scaling to obtain non-trivial limits when  $\max\{\th,\tilde\th\} < 1$. 
The first lemma indicates that it takes on the order of $n^2$ steps for the ERW to cross an interval of length $n$. 

\begin{lem}\label{lem:Tn}
 For any $k \in \Z$ let $T_k = \inf\{n\geq 0: X_n = k\}$. If $\max\{\th,\tilde\th\} < 1$, then there exist positive constants $\Cl[c]{ra}, \Cl[c]{ra2}>0$ such that
\begin{equation}\label{crossing}
 P\left( T_{\ell+n} - T_\ell \leq \frac{n^2}{L} \right) \leq \Cr{ra2} e^{-\Cr{ra} \sqrt{L}}
\quad \text{and} \quad 
 P\left( T_{-\ell-n} - T_{-\ell} \leq \frac{n^2}{L} \right) \leq \Cr{ra2} e^{-\Cr{ra} \sqrt{L}}, 
\end{equation}
for all integers $\ell\geq 0$, $n\geq 1$, and $L\in(0,\infty)$.
\end{lem}
\begin{proof}
  We shall show how to get the first inequality in \eqref{crossing},
  the proof of the second one being similar.  Since
  $T_\ell = \ell + 2 \sum_{y\leq \ell} \mathcal{D}_y^{(\ell,0)}$, it
  follows that
  $T_{\ell+n} - T_\ell \geq n + 2 \sum_{y=\ell}^n
  \mathcal{D}_y^{(\ell+n,0)}$.
  Moreover, the connection between $\mathcal{D}_y^{(\ell+n,0)}$ and BLP
  $\hat{V}$ (see the second to the last line of
    Table~\ref{table:ED_BLP} in Section \ref{connection}) implies
  that
\[
 P\left( T_{\ell+n} - T_\ell \leq \frac{n^2}{L} \right) \leq P\left( \sum_{i=0}^n \hat{V}_i \leq \frac{n^2}{2L} \, \biggl| \, \hat{V}_0 = 0 \right). 
 \]
We claim that it is sufficient to show that there exists an $L_0<\infty$ and $n_0\geq 1$ such that 
\begin{equation}\label{hVsum}
 P\left( \sum_{i=0}^n \hat{V}_i \leq \frac{n^2}{2L} \, \biggl| \, \hat{V}_0 = 0 \right) \leq e^{-\Cr{ra} \sqrt{L}}, \quad \forall L\geq L_0, \text{ and } n\geq \sqrt{L}n_0. 
\end{equation}
Indeed, since the statement of the Lemma holds trivially when $n < L$, this would then imply that $P( T_{\ell+n} - T_\ell \leq \frac{n^2}{L} ) \leq  e^{-\Cr{ra} \sqrt{L} }$ for $n\geq \max\{L_0,n_0^2\}$ and $L\geq L_0$, and then by choosing $\Cr{ra2}$ sufficiently large (depending on $L_0$ and $n_0$) the bound in \eqref{crossing} holds for all $n\geq 1$ and $0<L\leq n$. 
If $\theta \in (0,1)$, the proof of \eqref{hVsum} is the same as that of \cite[Lemma 3.2]{dkSLRERW}, where we use the first inequality of
Theorem~\ref{thm:BLPtail} (ii) instead of \cite[(2.2)]{dkSLRERW}.
On the other hand, if $\theta = s_{\hat{V}} \leq 0$ then this argument no longer works. We will thus handle the case $\theta\leq 0$ by coupling $\hat{V}$ with a ``smaller'' BLP $\hat{V}^{h,\e}$ which has parameter $\theta_{h,\e}  \in (0,1)$.

\textit{Reduction to the case $\theta\in(0,1)$ by coupling.}  Suppose
that $\w = \{\w_x(j)\}$ is the deterministic cookie environment with
periodic cookie stacks $\w_x = (p_1,p_2,\ldots,p_N,p_1,p_2,\ldots)$ as
given in Assumption \ref{asm:per}.  Fix an $h$ with
$0<h<\min_{i\leq N} (1-p_i)$ and for any $\e \in (0,1)$
let $\{G^\e_x\}_{x\in \Z}$ be an i.i.d.\ sequence of Geometric($\e$)
random variables; that is, $\P( G_x^\e = k) = (1-\e)^k \e$ for
$k\geq 0$.
Then, let $\w^{h,\e} = \{\w_x^{h,\e}(j)\}_{x\in\Z,\, j\geq 1}$ be a
random cookie environment constructed as follows.
\[
 \w^{h,\e}_x(j) = \begin{cases} \w_x(j) + h & \text{if } j\leq G_x^\e \\ \w_x(j) & \text{if } j > G_x^\e,  \end{cases}
\qquad x \in \Z, \, j \geq 1. 
\]
Since the above construction couples the cookie environments $\w$ and $\w^{h,\e}$ in such a way that $\w_x(j) \leq \w_x^{h,\e}(j)$, it follows that we can couple BLPs $\hat{V}$ and $\hat{V}^{h,\e}$ the the corresponding cookie environments so that $\hat{V}_i \geq \hat{V}^{h,\e}_i$ for all $i\geq 0$. 
The random cookie environment $\w^{h,\e}$ fits into the framework of the Markovian cookie stacks considered in \cite{kpERWMC}; in particular, by \cite[Theorem 2.7]{kpERWMC} the tail asymptotics in Theorem \ref{thm:BLPtail} hold for $Z=\hat{V}^{h,\e}$ with parameter $s_Z = \theta_{h,\e}$ that can be calculated using \cite[Lemma 5.1]{kpERWMC} to be $\theta_{h,\e} = \theta + \frac{4h(1-\e)}{\nu\e}$. 
Therefore, if $\theta \leq 0$ we may choose an $\e\in (0,1)$ so that $\theta_{h,\e} \in (0,1)$. Then the coupling of $\hat{V}$ and $\hat{V}^{h,\e}$ together with the argument above for $\theta \in (0,1)$ applied to the BLP $\hat{V}^{h,\e}$ will ensure that 
\[
 P\left( \sum_{i=0}^n \hat{V}_i \leq \frac{n^2}{2L} \, \biggl| \, \hat{V}_0 = 0 \right)
\le P\left( \sum_{i=0}^n \hat{V}^{h,\e}_i \leq \frac{n^2}{2L} \, \biggl| \, \hat{V}^{h,\e}_0 = 0 \right) \leq e^{-\Cr{ra} \sqrt{L}}.
\]

\end{proof}

Lemma \ref{lem:Tn} easily implies the following control for the range of the EWR. 
\begin{cor}\label{range}
 If $\max\{\th,\tilde\th\} < 1$, then 
\[
 P\left( \sup_{k\leq n} |X_k| > K \sqrt{n} \right)\le 2 \Cr{ra2} e^{-\Cr{ra}K}>0, 
\forall n\geq 1, \, K> 0,  
\]
where $\Cr{ra},\Cr{ra2}$ are the constants from Lemma \ref{lem:Tn}.
\end{cor}

Another consequence of Lemma \ref{lem:Tn} is the following process-level tightness estimate for the running minimum and maximum of the ERW. 
\begin{cor}\label{cor:MnIntight} 
Let $M_n = \sup_{k\leq n} X_k$ and $I_n = \inf_{k\leq n} X_k$ be the running maximum and minimum, respectively, of the excited random walk. 
If $\max\{ \th, \tilde\th \} < 1$, then for any $\e>0$ and $t<\infty$, 
\[
 \lim_{\d \ra 0} \limsup_{n\ra\infty} P\left( \sup_{\substack{k,\ell\leq nt\\ |k-\ell|\leq n\d}} |M_k - M_\ell | \geq 2\e \sqrt{n} \right) = 
\lim_{\d \ra 0} \limsup_{n\ra\infty} P\left( \sup_{\substack{k,\ell\leq nt\\ |k-\ell|\leq n\d}} |I_k - I_\ell | \geq 2\e \sqrt{n} \right) =
0.
\]
\end{cor}
\begin{proof}
 It's enough to prove the limit for the running maximum process since the proof is the same for the running minimum. 
If the running maximum increases at least $2 \e \sqrt{n}$ over some time interval less than $\d n$, then it follows that some interval of the form $[(m-1)\fl{\e \sqrt{n}}, m \fl{\e \sqrt{n}}]$ is crossed in less than $\d n$ steps. 
Therefore, 
\begin{align*}
  P\left( \sup_{\substack{k,\ell\leq nt\\ |k-\ell|\leq n\d}} |M_k - M_\ell | \geq \e \sqrt{n} \right)
&\leq P\left( M_n \geq \frac{\e \sqrt{n}}{\d} \right) + \sum_{m=1}^{\fl{1/\d}} P\left( T_{m \fl{\e \sqrt{n}}} - T_{(m-1)\fl{\e\sqrt{n}}} \leq \d n \right), 
\end{align*}
and the conclusion of Corollary \ref{cor:MnIntight} follows easily from this and Lemma \ref{lem:Tn}. 
\end{proof}


\begin{lem}\label{lt}
Let $\mathcal{L}(n;x) = \sum_{k=0}^{n-1} \ind{X_k = x}$ be the number of visits to $x$ prior to time $n$.
If $\max\{\th,\tilde\th\} < 1$, then 
\[
 \lim_{K\ra\infty} \limsup_{n\ra\infty} P\left( \sup_{x\in \Z} \mathcal{L}(n;x) > K \sqrt{n} \right) = 0. 
\]
\end{lem}

\begin{proof}
By symmetry it is enough to show that 
\[
 \lim_{K\ra\infty} \limsup_{n\ra\infty} P\left( \sup_{x\geq 0} \mathcal{L}(n;x) > K\sqrt{n} \right) = 0.  
\]
If some site $x\geq 0$ is visited more than $K\sqrt{n}$ times by time $n$, then during the first $\fl{\sqrt{K n}}$ excursions to the right of the origin either 
\begin{itemize}
 \item the walk takes at most $n$ steps to complete these excursions to the right,
 \item or some point to the right of the origin is visited at least $K \sqrt{n}$ times during these excursions. 
\end{itemize}
For any $m\geq 1$, the sum $2\sum_{x\geq 0} \mathcal{E}_x^{(0,m)}$ gives the total amount of time taken by the excursions to the right of the origin during the first $m$ excursions from the origin (to the right or left). 
Also, during the first $m$ excursions from the origin the total time spent at a site $x\geq 0$ is equal to  $\mathcal{D}_x^{(0,m)} + \mathcal{E}_x^{(0,m)} = \mathcal{E}_{x-1}^{(0,m)} + \mathcal{E}_x^{(0,m)}$. Therefore, 
\begin{align}
& P\left( \sup_{x\geq 0} \mathcal{L}(n;x) > K\sqrt{n} \right) \nonumber  \\
&\qquad \leq P\left( 2\sum_{x\geq 0} \mathcal{E}_x^{(0,m)} \leq n \, \biggl| \, \mathcal{E}_0^{(0,m)} = \fl{\sqrt{Kn}} \right) + P\left( \sup_{x \geq 0} \mathcal{E}_x^{(0,m)} \geq \frac{ K \sqrt{n} }{2} \, \biggl| \, \mathcal{E}_0^{(0,m)} = \fl{\sqrt{Kn}} \right) \nonumber \\
&\qquad = P\left( 2 \sum_{i=0}^{\s_0^U} U_i \leq n \, \biggl| \, U_0 = \fl{\sqrt{Kn}} \right) + P\left( \sup_{0\leq i < \s_0^U} U_i \geq \frac{K\sqrt{n}}{2} \, \biggl| \, U_0 = \fl{\sqrt{Kn}} \right), \label{ltblp}
\end{align}
where in the last line we used the connection with the BLP detailed in Section \ref{connection}. 
For the first probability in \eqref{ltblp}, note that the diffusion approximation in Theorem \ref{thm:da} implies that for any $\e>0$, 
\begin{align*}
 \limsup_{n\ra\infty} P\left( 2 \sum_{i=0}^{\s_0^U} U_i \leq n \, \biggl| \, U_0 = \fl{K\sqrt{n}} \right) 
&\leq \lim_{n\ra\infty} P\left( 2 \sum_{i=0}^{\s_{\e \sqrt{n}}^U} U_i \leq n \, \biggl| \, U_0 = \fl{K\sqrt{n}} \right) \\
&= P\left( 2 \int_0^{\s_\e^Y} Y_s \, ds \leq \frac{1}{K^2} \, \biggl| \, Y_0 = 1 \right), 
\end{align*}
where $Y(t)$ is solves the SDE
$dY(t) = \rho \, dt + \sqrt{\nu Y(t)} \, dB(t)$.  The last probability
vanishes when we let $K\to\infty$. 
The second probability in \eqref{ltblp} can be bounded
  uniformly in $n$ by a quantity that vanishes as $K\ra\infty$. This
  is asserted by the following lemma.
\begin{lem}
  \label{5.5}
If $\th<1$, then for every $\epsilon>0$ there is a constant
  $\Cl[c]{km}=\Cr{km}(\epsilon)$ such that
\[
 P\left( \sup_{i< \s_0^U} U_i > \Cr{km} n \, \biggl| \, U_0=n \right) < \e, \quad \forall n\in \N. 
\]
\end{lem}
\noindent The proof of this lemma is the same as that of \cite[(5.5)]{kmLLCRW}
and uses the analogs of \cite[Lemmas 6.3,\ 6.4]{kpERWMC} for the
process $U$ in place of \cite[Lemmas 5.1,\ 5.3]{kmLLCRW} respectively.
\end{proof}

\subsection{Control of rarely visited sites}
Lemmas \ref{lem:Tn} and \ref{lt} together imply that there are of the order $\sqrt{n}$ sites with local time of the order $\sqrt{n}$. However, there may be some sites in the range that have been visited far fewer than $\sqrt{n}$ times. 
The lemmas in this section give control on how often one of the BLPs can be below some fixed level. 
Since the BLPs are related to the local times of directed edges, these then give control on the number of directed edges which have been traversed a relatively small number of times. 
 
\begin{lem}\label{lem:BLPsmall1}
  Let $Z$ be one of the BLPs $U$, $\hat{U}$, $V$, $\hat{V}$ with
  $s_Z>0$, where $s_Z$ is defined in
    Theorem~\ref{thm:BLPtail}. For every $\gamma \in (0,1/2)$ and
  $\epsilon>0$ there exist positive constants
  $\Cl[c]{sm1}=\Cr{sm1}(\gamma,\epsilon)$ and
  $\Cl[c]{sm2}=\Cr{sm2}(\gamma,\epsilon)$ such that
\[
P\left( \sum_{i=1}^{\s^Z_0} \ind{Z_{i-1} < n^{\gamma}} > \e \sqrt{n} \, \biggl| \, Z_0 = m \right) \le  \Cr{sm1}e^{-\Cr{sm2} n^{1/2-\gamma}},  \qquad \forall m\geq 0, \, n\geq 1. 
\]
\end{lem}

\begin{proof}
  By the Markov property, it is enough to prove the statement of the
  lemma 
for $m<n^\gamma$.
Let
  $k=\min\{j\in\N:\,2^j\ge n^\gamma\}$. Define
  $I_j=[2^{k-j},2^{k-j+1}),\ j\in\N$, and events
  \[ A_j
=\left\{\sum_{i=1}^{\sigma_0^Z}\ind{Z_{i-1}\in I_j}>\frac{\epsilon j(2^{j-1}|I_j|)^{1/(2\gamma)}}{2^{j+1}}\right\}.
\]
  Since $2^{j-1}|I_j|=2^{k-1}< n^\gamma$ and
  $\sum_{j=1}^kj2^{-(j+1)}<1$, we get the inclusion 
  \[\left\{\sum_{i=1}^{\s^Z_0} \ind{Z_{i-1} < n^{\gamma}} > \e \sqrt{n}
  \right\}\subseteq\bigcup_{j=1}^k A_j.\]
  Thus, it is enough to show that
  $\sum_{j=1}^k P(A_j\,|\,Z_0=m)\le
  \Cr{sm1}e^{-\Cr{sm2}n^{1/2-\gamma}}$.
  To estimate the probabilities of sets $A_j$ we shall need the
  following proposition, which is an adaptation of \cite[Proposition
  6.1]{kmLLCRW}.
  \begin{prop}\label{6.1}
    Let $Z$ be one of the BLPs $U$, $\hat{U}$, $V$, $\hat{V}$  with
    $s_Z>0$. Then there is a positive constant $\Cl[c]{pr}$ such that
    for all $n,x\in\N$, and $m\ge 0$
    \[P\left(\sum_{i=1}^{\sigma_0^Z}\ind{Z_{i-1}\in[x,2x)}>2xn\,\Big|\,Z_0=m\right)\le
    e^{-\Cr{pr}n}.\]
  \end{prop}
\begin{proof}[Proof of Proposition~\ref{6.1}]
 The proof for all four processes is identical to that
    of Proposition 6.1 in \cite{kmLLCRW}. For reader's convenience
    we note that our process $\hat{V}$ corresponds to the process $V$
    in \cite{kmLLCRW} and $s_Z$ corresponds to $\delta$ in
    \cite{kmLLCRW}. In the proof of Proposition 6.1 we only need to
    replace the references to Lemma 3.1 and Lemma 5.3 of
    \cite{kmLLCRW} with the references to Theorem~\ref{thm:da} above
    and Lemma~6.4 of \cite{kpERWMC} respectively. The only place which
    uses the inequality $s_Z>0$ is Corollary 5.5 of
    \cite{kmLLCRW}. Since this corollary depends only on \cite[Lemma
    5.3]{kmLLCRW}, which is fully replaced by Lemma~6.4 of
    \cite{kpERWMC} in our setting, the proof goes through without any
    changes. 
\end{proof}
Applying Proposition~\ref{6.1} to each term of the sum we get
\begin{align*}
 \sum_{j=1}^k P(A_j \, | \, Z_0 = m) &\leq \sum_{j=1}^k \exp\left( -\Cr{pr}\left\lfloor \frac{\e j 2^{(j-1)/(2\gamma)} |I_j|^{1/(2\gamma)-1} }{2^{j+2}} \right\rfloor \right) \\
&=  \sum_{j=1}^k \exp\left( -\Cr{pr}\left\lfloor \e j 2^{k(\frac{1}{2\gamma}-1)-\frac{1}{2\gamma}-2} \right\rfloor \right)
\leq  \sum_{j=1}^\infty \Cl[c]{prnb} \exp\left(-\Cl[c]{prn}\epsilon j n^{1/2-\gamma}\right).
\end{align*}
This immediately implies the statement of the lemma.
\end{proof}
\begin{lem}\label{lem:BLPsmall2}
 Let $Z$ be either the BLP $\hat{U}$ or $\hat{V}$ and let $\max\{\th,\tilde\th\} < 1$. 
For every
  $\gamma \in (0,1/2)$, $\epsilon>0$ and $\beta>1$ there is a
  $\Cl[c]{sm3}=\Cr{sm3}(\gamma,\beta,\epsilon)$ such
  that for any $K>1$ and all sufficiently large $n$
\[
 \max_{m\ge 0} P\left( \sum_{i\le K\sqrt{n}} \ind{ Z_{i-1} < n^{\gamma} } > \e \sqrt{n} \, \biggl| \, Z_0 = m \right) \le \frac{\Cr{sm3}}{n^\beta}. 
\]
\end{lem}

\begin{proof}[Proof of Lemma \ref{lem:BLPsmall2}]
  We will give the proof only for the process $\hat{V}$
    as the proof for $\hat{U}$ is similar.  First of all, by the
    monotonicity of the BLPs with respect to their initial conditions,
    the maximum is attained at $m=0$. Thus, we shall set
    $m=0$. Secondly, we note that it is sufficient to consider only
    the case $\theta=s_{\hat{V}}\in(0,1)$. Indeed, suppose that
    $\theta\le 0$. Just as in the proof of Lemma~\ref{lem:Tn}, we can
    invoke \cite[Lemma 5.1]{kpERWMC} to construct a ``smaller'' BLP process
    $\hat{V}^{h,\epsilon}$ such that
    $\theta_{h,\epsilon}=s_{\hat{V}^{h,\epsilon}}  \in(0,1)$ and
    $P(\forall i\ge 0\ \ \hat{V}^{h,\epsilon}_i  \le \hat{V}_i)=1$.  Then
    \[P\left( \sum_{i\le K\sqrt{n}} \ind{ \hat{V}_{i-1} < n^{\gamma} }
      > \e \sqrt{n} \, \biggl| \, \hat{V}_0 = 0 \right)\le P\left(
      \sum_{i\le K\sqrt{n}} \ind{ \hat{V}^{h,\epsilon}_{i-1}  < n^{\gamma} } > \e
      \sqrt{n} \, \biggl| \,  \hat{V}^{h,\epsilon}_0  = 0 \right),\]
    and it is sufficient to show that the last probability does not exceed $\Cr{sm3}/n^\beta$. Thus,
    without loss of generality we shall also assume that $\theta\in(0,1)$.

{\em Step 1.} We will first estimate the
  number of times $\hat{V}_i = 0$ for $i\leq K\sqrt{n}$. Let
  $\sigma^{\hat{V}}_{0,0}=0$ and denote by $\sigma^{\hat{V}}_{0,i} = \inf \{j> \s_{0,i-1}^{\hat{V}} : \, \hat{V}_j = 0 \}$
  the $i$-th hitting time of $0$. If $\sigma^{\hat{V}}_{0,i-1}=\infty$ for some $i\in\N$ then we set $\sigma^{\hat{V}}_{0,j}=\infty$ for all $j\ge i$ and $\sigma^{\hat{V}}_{0,i}-\sigma^{\hat{V}}_{0,i-1}=\infty$. Define
  $\alpha=(1-\theta)/4$. Then 
\begin{multline*}
P\left( \s_{0,\fl{n^{1/2-\a}}}^{\hat{V}} \leq K \sqrt{n} \right)
\le \prod_{i\le n^{1/2-\alpha}} P(\sigma^{\hat{V}}_{0,i}-\sigma^{\hat{V}}_{0,i-1} \le K\sqrt{n}\,|\,\hat{V}_0=0)\\\overset{\text{Th.~\ref{thm:BLPtail}}}{\le}
\left(1-\frac{C^{\hat{V}}_1(0)}{2\fl{K\sqrt{n}}^{\theta}}\right)^{\fl{n^{1/2-\alpha}}}\le
  \exp\left(-\Cl[c]{sm6}n^{(1-\theta-2\alpha)/2}\right),
\end{multline*}
for some positive constant $\Cr{sm6}=\Cr{sm6}(\theta,K)$. 

  {\em Step 2.} Now that we know that the number of regenerations can
  not be too large, we can add up the time spent below level $n^{\gamma}$
  for up to $\fl{n^{1/2-\alpha}}$ regenerations.  Let $t_i=\sum_{j=\sigma^{\hat{V}}_{0,i-1}}^{\sigma^{\hat{V}}_{0,i}-1}\ind{\hat{V}_j<n^{\gamma}}$, $i\in\N$. Then
\begin{multline*}
P\left( \sum_{i\le K\sqrt{n}} \ind{\hat{V}_{i-1} < n^{\gamma}} > \epsilon\sqrt{n} \,\Big|\,\hat{V}_0=0\right)\\
\le P\left( \s_{0,\fl{n^{1/2-\a}}}^{\hat{V}} \leq K \sqrt{n} \right)
    +P\left(\sum_{i\le n^{1/2-\alpha}}t_i>\epsilon\sqrt{n} \,\Big|\,\hat{V}_0=0\right).
\end{multline*}
Note that the $t_i$'s are independent and identically
distributed. Lemma~\ref{lem:BLPsmall1} (for $Z=\hat{V}$) implies that
the tail of each $t_i$ decays faster than any power of $n$, so every
moment of $t_i$ is finite.  Using Step 1, the Markov inequality with
power $\ell > \beta/\alpha$, and then Jensen's inequality we get
\begin{align*}
 P\left( \sum_{i\le K\sqrt{n}} \ind{\hat{V}_{i-1} < n^{\gamma} } > \epsilon\sqrt{n} \,\Big|\,\hat{V}_0=0\right)
&\le P\left( \s_{0,\fl{n^{1/2-\a}}}^{\hat{V}} \leq K \sqrt{n} \right) +\frac{ E\left[\left(\sum_{i\le n^{1/2-\alpha}}t_i\right)^\ell\,\Big|\,\hat{V}_0=0\right]}{\epsilon^\ell n^{\ell/2}} \\
&\le \exp \left( -\Cr{sm6}n^{(1-\theta-2\alpha)/2} \right) +\frac{ n^{(1/2-\a)\ell} E\left[t_1^\ell\,\Big|\,\hat{V}_0=0\right]}{\epsilon^\ell n^{\ell/2}} \le \frac{\Cr{sm3}}{n^\beta}, 
\end{align*}
for $n$ large enough as claimed. 
\end{proof}

\begin{cor}\label{cor:ltsmall} Let $\max\{\theta,\hat{\theta}\}<1$. For every
  $\gamma\in(0,1/2)$, $\epsilon>0$ and $\beta>1$ there is a
  $\Cl[c]{sm7}=\Cr{sm7}(\gamma,\beta,\epsilon)$ such
  that for all sufficiently large $n$
\[
 P\left( \sum_{y \in \Z} \ind{1 \leq \mathcal{L}(n;y) < n^{\gamma}} > \e \sqrt{n} \right)\le\frac{\Cr{sm7}}{n^\beta}.
\]
\end{cor}
\begin{proof}
Suppose that after $n$ steps the ERW is at $X_n=x$ and has previously visited that site $m$ times; 
that is, $\l_{x,m}=n$. In this case, the local times at sites can be expressed using the directed edge local times $\mathcal{E}^{(x,m)}_y$ and $\mathcal{D}^{(x,m)}_y$ as given in \eqref{LED}. 
In particular, if $0\leq x\leq n$ and $m\geq 0$ this implies that 
\begin{align*}
& P\left( \sum_{y \in \Z} \ind{1 \leq \mathcal{L}(n;y) < n^{\gamma}} > \e \sqrt{n}, \, \l_{x,m} = n \right) \\
&\leq P\left( \sum_{y\leq 0} \ind{1 \leq \mathcal{D}^{(x,m)}_y < n^\gamma}  + \sum_{0 < y \leq x} \ind{\mathcal{D}^{(x,m)}_y < n^\gamma} + \sum_{y\geq x} \ind{1 \leq \mathcal{E}^{(x,m)}_y < n^\gamma} > \e \sqrt{n}, \, \l_{x,m} = n \right) \\
&\leq \max_{m'\geq 0} P\left( \sum_{i=0}^{\s_0^V-1} \ind{V_i < n^\gamma} > \frac{\e \sqrt{n}}{3} \, \biggl| \, V_0 = m' \right) 
+ \max_{m'\geq 0} P\left( \sum_{i\leq n} \ind{\hat{V}_i < n^\gamma} > \frac{\e \sqrt{n}}{3} \, \biggl| \, \hat{V}_0 = m' \right) \\
&\qquad + \max_{m'\geq 0} P\left( \sum_{i=0}^{\s_0^U-1} \ind{U_i < n^\gamma} > \frac{\e \sqrt{n}}{3} \, \biggl| \, U_0 = m' \right). 
\end{align*}
Note that this upper bound is uniform over all $0\leq x \leq n$ and $m\geq 0$, and that Lemmas \ref{lem:BLPsmall1} and \ref{lem:BLPsmall2} imply that this is bounded above by $\Cr{sm3}(\gamma,\b+2,\e/3) n^{-\b-2}$ for all $n$ large enough. A similar upper bound holds uniformly for $-n\leq x<0$ and $m\geq 0$. Finally, since exactly one of the events $\{\l_{x,m} = n\}$ occurs for $|x|\leq n$ and $m<n$, we have that 
\[
  P\left( \sum_{y \in \Z} \ind{1 \leq \mathcal{L}(n;y) < n^{\gamma}} > \e \sqrt{n} \right)
= \sum_{|x|\leq n, \, m<n} P\left( \sum_{y \in \Z} \ind{1 \leq \mathcal{L}(n;y) < n^{\gamma}} > \e \sqrt{n}, \, \l_{x,m} = n \right) 
\leq \frac{\Cr{sm7}}{n^\b},
\]
for all $n$ large enough. 
\end{proof}

\section{Functional limit laws: Non-boundary cases}\label{sec:pbm}
In this section we give the proof of Theorem \ref{thm:pbm}: convergence to perturbed Brownian motion for ERW when $\max\{\th,\tilde\th\} < 1$. 
Our proof mimics in some ways the proof in \cite{dkSLRERW} for the case of ERW with boundedly many cookies per site in that we decompose the walk into a martingale plus a term that is approximately equal to a linear combination of the running maximum and minimum of the walk. However, the decomposition in this paper is different and less transparent than in \cite{dkSLRERW} (as will be evident from the proof below, this is due to the fact that $\theta+\tilde\theta$ is not necessarily equal to zero; in contrast, for ERW with boundedly many cookies per site we have $\theta = \d$ and $\tilde\th = -\d$). 
Moreover, the control of the non-martingale part of the decomposition is much more difficult than in \cite{dkSLRERW} due to the fact that there may be infinitely many cookies at a site with non-zero drift. 

\textbf{Step 1: Control of martingale term.}
We begin by noting that if $\mathcal{F}_n= \s(X_k, \, k\leq n)$ then $E[ X_{n+1}-X_n | \mathcal{F}_n] = 2\w_{X_n}\left(\mathcal{L}(n+1;X_n)\right) - 1$. Therefore, if 
\begin{equation}\label{Cndef}
 C_n = \sum_{k=0}^{n-1} (2\w_{X_k}(\mathcal{L}(k+1;X_k))-1) = \sum_{y \in \Z} \sum_{j=1}^{\mathcal{L}(n;y)} (2\w_y(j)-1),
\end{equation}
it follows that $B_n = X_n - C_n$ is a martingale with respect to the filtration $\mathcal{F}_n$. 
The first step in the proof of Theorem \ref{thm:pbm} is the following control of this martingale. 

\begin{lem}\label{lem:BnBM}
Let $B_n = X_n - C_n$, with $C_n$ defined in \eqref{Cndef}. Then, 
\begin{equation}\label{Bnlimit}
 \left\{\frac{B_{tn}}{\sqrt{(\nu/2) n}}\right\}_{t\geq 0} \xRightarrow[n\ra\infty]{J_1} \{ \mathfrak{B}(t) \}_{t\geq 0}, 
\end{equation}
where $\mathfrak{B}(t)$ is a standard Brownian motion and $\nu$ is the parameter given in \eqref{rhonuforms}. 
\end{lem}

\begin{proof}
Since $B_n$ is a martingale with bounded steps we need only to check the convergence of the quadratic variation of the process \cite[Theorem 18.2]{bCOPM}.
In particular, it is sufficient to show that 
\begin{equation} \label{Bnqv}
 \lim_{n\ra\infty} \frac{1}{n} \sum_{k\leq n} E\left[ (B_k-B_{k-1})^2 \, | \, \mathcal{F}_{k-1} \right] = \frac{\nu}{2}, \qquad P\text{-a.s.} 
\end{equation}
Note first of all that 
\begin{align*}
 \sum_{k=1}^{n} E\left[ (B_k-B_{k-1})^2 \, | \, \mathcal{F}_{k-1} \right] 
&= \sum_{k=1}^{n} \left\{ E\left[( X_k-X_{k-1})^2 \,|\, \mathcal{F}_{k-1} \right] - E\left[ X_k-X_{k-1} \,|\,\mathcal{F}_{k-1}  \right]^2 \right\} \\
&= \sum_{k=1}^{n} \left\{ 1 - \left( 2 \w_{X_{k-1}}(\mathcal{L}(k;X_{k-1})) - 1 \right)^2 \right\} \\
&= n - \sum_{y\in \Z} \sum_{j=1}^{\mathcal{L}(n;j)} (2\w_y(j) - 1)^2. 
\end{align*}
Therefore, \eqref{Bnqv} is equivalent to 
\begin{equation}\label{Bnqv2}
 \lim_{n\ra\infty} \frac{1}{n} \sum_{y \in\Z}\sum_{j=1}^{\mathcal{L}(n;y)} (2\w_y(j)-1)^2 = \frac{1}{N} \sum_{i=1}^N (2p_i-1)^2 = 1-\frac{\nu}{2}, \qquad P\text{-a.s.},
\end{equation}
where the last equality follows easily from the formula for $\nu$ in \eqref{rhonuforms}.
To prove the first equality in \eqref{Bnqv2}, we first note that since $\{\w_y(j)\}_{j\geq 1}$ is a deterministic periodic sequence then 
\[
 \left| \frac{1}{m} \sum_{j=1}^m (2\w_y(j) - 1)^2 - \frac{1}{N}\sum_{i=1}^N (2p_i-1)^2 \right| \leq \frac{C}{m}, \quad \forall m\geq 1, 
\]
for some $C>0$. Since $\sum_y \mathcal{L}(n;y) = n$, it follows that 
\begin{align*}
& \left|\frac{1}{n} \sum_{y \in\Z}\sum_{j=1}^{\mathcal{L}(n;y)} (2\w_y(j)-1)^2 - \frac{1}{N} \sum_{i=1}^N (2 p_i-1)^2 \right| \\
&=  \left| \frac{1}{n} \sum_{y: \mathcal{L}(n;y) \geq 1} \mathcal{L}(n;y) \left\{ \frac{1}{\mathcal{L}(n;y)} \sum_{j=1}^{\mathcal{L}(n;y)}  (2\w_y(j)-1)^2 - \frac{1}{N} \sum_{i=1}^N (2 p_i-1)^2 \right\} \right| \\
&\leq \frac{1}{n} \sum_{y: \mathcal{L}(n;y) \geq 1} \frac{C}{\mathcal{L}(n;y)}\le
\frac{C}{n^\gamma} +  \frac{C}{n} \sum_{y\in\Z} \ind{1\leq \mathcal{L}(n;y) \leq n^\gamma}, 
\end{align*}
for any $\gamma \in (0,1/2)$.  It follows from Corollary
\ref{cor:ltsmall} and the Borel-Cantelli Lemma that the sum in the
last line is $o(\sqrt{n})$, $P$-a.s., and from this the limit in
\eqref{Bnqv2} follows.
\end{proof}

\textbf{Step 2: Control of accumulated drift.}
Having proved Lemma \ref{lem:BnBM} we now need to control the term $C_n$ defined in \eqref{Cndef} which records the total drift which the ERW has accumulated from the cookie environment. The following lemma shows that $C_n$ is approximated by a fixed linear combination of the distance of the random walk from its running maximum and running minimum. 
Recall that $M_n$ and $I_n$ were defined in Corollary \ref{cor:MnIntight} as the running maximum and minimum, respectively, of the ERW. 
\begin{lem}\label{lem:Cnapprox}
 For any $t\geq 0$ and any $\e>0$, 
\[
\lim_{n\ra\infty} P\left( \sup_{k\leq nt} \left| C_k - \rho \left(M_k - X_k \right) - \tilde\rho\left(I_k - X_k \right) \right| \geq \e \sqrt{n} \right)  = 0, 
\]
where $\rho$ and $\tilde\rho$ are the parameters defined in \eqref{BLPrhonu}. 
\end{lem}

\begin{proof}
For the proof of Lemma \ref{lem:Cnapprox}, we need to show that total drift contained in the used cookies to the right (resp. left) of $X_k$ is approximately $\rho$ (resp. $-\tilde\rho$) times the number of sites visited to the right (resp. left) of $X_k$. 
That is, if we decompose $C_n$ as $C_n = C_n^-+C_n^0 + C_n^+$, where 
\[
 C_n^+ = \sum_{y > X_n} \sum_{j=1}^{\mathcal{L}(n;y)} (2\w_y(j)-1), \quad C_n^- = \sum_{y < X_n} \sum_{j=1}^{\mathcal{L}(n;y)} (2\w_y(j)-1), 
\quad \text{and} \quad C_n^0 = \sum_{j=1}^{\mathcal{L}(n;X_n)} (2\w_0(j)-1), 
\]
then since  $C_n^0 \in \{ \sum_{j=1}^\ell (2p_j-1) , \, \ell=1,2,\ldots,N\}$ is a bounded random variable it is enough to show
\begin{equation}\label{Cplus}
 \lim_{n\ra\infty} P\left( \sup_{k\leq nt} \left| C_k^+ - \rho(M_k - X_k) \right| \geq \e \sqrt{n} \right) = 0, \qquad \forall \e>0, \, t<\infty, 
\end{equation}
and
\begin{equation}\label{Cminus}
 \lim_{n\ra\infty} P\left( \sup_{k\leq nt} \left| C_k^- - \tilde\rho(I_k - X_k) \right| \geq \e \sqrt{n} \right) = 0, \qquad \forall \e>0, \, t<\infty.
\end{equation}

The proofs of \eqref{Cplus} and \eqref{Cminus} are similar, and thus we will only give the proof of \eqref{Cplus}. We will prove this using properties of the BLP and the connection with the random walk given in Section \ref{connection}. 
To this end, for any $x,y\in \Z$ and $m\geq 0$ let 
\begin{equation}\label{Dyxm}
 \Delta_y^{(x,m)} = \sum_{j=1}^{\mathcal{L}(\l_{x,m};y)} (2\w_y(j) -1), 
\end{equation}
be the total drift contained in the cookies used at site $y$ prior to time $\l_{x,m}$ (recall this is the time of the $(m+1)$-st visit to site $x$). 
With this notation we have that $C_k^+ = \sum_{y>x} \Delta_y^{(x,m)}$ on the event $\{ X_k = x, \, \mathcal{L}(k;x) = m\} = \{\l_{x,m} = k\}$.
Note that since $\Delta_y^{(x,m)} = 0$ for sites $y$ that have not been visited, we can restrict the sum in this representation of $C_k^+$ to $y\leq \mathcal{M}^{(x,m)} := M_{\l_{x,m}}$. 
Therefore, with this notation we have that 
\begin{equation}\label{CplusDelta}
 C_k^+ - \rho(M_k-X_k) = \sum_{y=x+1}^{\mathcal{M}^{(x,m)}} \left( \Delta_y^{(x,m)} - \rho \right), \qquad \text{on the event } \{\l_{x,m} = k\}. 
\end{equation}
We will use the representation in \eqref{CplusDelta} to prove \eqref{Cplus} by showing that the terms inside the sum on the right in \eqref{CplusDelta} are nearly equal to a martingale difference sequence. 
To this end, for $x\in \Z$ and $m\geq 0$ fixed let
\[
 \rho^{(x,m)}_y = E\left[ \Delta_y^{(x,m)} \, | \, \mathcal{G}_{y-1}^{(x,m)} \right], 
\qquad \text{where } 
\mathcal{G}^{(x,m)}_z = \s\left( \mathcal{E}_y^{(x,m)}: \,y \leq z \right). 
\]
Then $\{\Delta_y^{(x,m)} - \rho_y^{(x,m)} \}_{y>x}$ is a martingale difference sequence with respect to the filtration $\{\mathcal{G}_y^{(x,m)}\}_{y>x}$, and the following lemma shows that this is not too far from the original sequence. 
\begin{lem}\label{rhoyxm}
 There exist constants $C,c>0$ such that $\left| \rho_y^{(x,m)} - \rho \right| \leq C \exp\left\{-c\, \mathcal{E}_{y-1}^{(x,m)} \right\}$ for all $x \in \Z$, $m\geq 0$ and $y> x$.   
\end{lem}

\begin{proof}
It follows from the definition of $\Delta_y^{(x,m)}$ in \eqref{Dyxm}, the connection between local time a sites and the directed-edge local times in \eqref{LED}, and the connection of $\{\mathcal{E}^{(x,m)}_{y}\}_{y\geq x}$ with the BLPs $U$ and $\hat{U}$ that
\[
 \rho^{(x,m)}_y = E\left[ \sum_{j=1}^{ \mathcal{E}^{(x,m)}_{y} + \mathcal{E}^{(x,m)}_{y-1} + \ind{x<y\leq 0}} (2\w_y(j)-1) \, \biggl| \, \mathcal{E}_{y-1}^{(x,m)} \right]
= \psi\left( \mathcal{E}_{y-1}^{(x,m)} + \ind{x<y\leq 0} \right), 
\]
where
$\psi(n) = E\left[ \sum_{j=1}^{U_1+U_0} (2 \w_1(j)-1) \, \bigl| \, U_0
  = n \right]$.
Therefore, we need only to show that $|\psi(n) - \rho| \leq C e^{-cn}$ for some constants $C,c>0$. 
To this end, first note that 
\begin{align}
 \psi(n) &= E\left[ \sum_{j=1}^{U_1+n} \left\{ 2( \w_1(j)- \xi_1(j))-2(1- \xi_1(j)) + 1 \right\} \biggl| \, U_0 = n \right] \nonumber \\
&= 2 E\left[ \sum_{j=1}^{U_1+n} ( \w_1(j)- \xi_1(j)) \biggl| \, U_0 = n \right]
- 2 E\left[ \sum_{j=1}^{U_1+n} (1-\xi_1(j)) \biggl| \, U_0 = n \right] + E[U_1 | \, U_0 = n] + n. \label{psinrep}
\end{align}
By the construction of the BLP $U$, $U_1 + n$ is the number of trials in the Bernoulli sequence $\{\xi_1(j)\}_{j\geq 1}$  until the $n$-th failure. Thus, the sum inside the second expectation in \eqref{psinrep} equals $n$ and the Optional Stopping theorem implies that the first expectation in \eqref{psinrep} is equal to zero. 
Thus, we have shown that $\psi(n) = E[U_1|\, U_0=n] - n$ and it follows from 
\eqref{rhoroc}
that $|\psi(n)-\rho| \leq C e^{-cn}$. 
\end{proof}

We are now ready to give the proof of \eqref{Cplus}.
For any $n\geq 1$ and $K,t<\infty$ let 
\[
 G_{n,K,t} = \left\{ \sup_{k\leq nt} |X_k| \leq K \sqrt{n} \right\} \cap \left\{ \sup_{y \in \Z} \mathcal{L}(\fl{nt}; y) \leq K\sqrt{n} \right\}
 \cap \bigcap_{\substack{|x|\leq K\sqrt{n} \\ m\leq K\sqrt{n}}} \left\{ \sum_{y=x+1}^{\mathcal{M}^{(x,m)}} \ind{\mathcal{E}_{y-1}^{(x,m)} < n^{1/4} } > \frac{\sqrt{n}}{K} \right\}. 
\]
Note that Corollary \ref{range} and Lemma \ref{lt} imply that the first two events on the right are typical events for $K$ sufficiently large. To see that the intersection of the other events is also typical, recall that for any fixed $(x,m)$ the directed edge local time process $\{\mathcal{E}^{(x,m)}_y\}_{y\geq x}$ is a Markov chain with transition probabilities given by the BLP $U$ and $\hat{U}$. In particular, for any fixed $|x| \leq K\sqrt{n}$ 
\begin{align*}
\sup_{m\geq 0} P\left( \sum_{y=x+1}^{\mathcal{M}^{(x,m)}} \ind{\mathcal{E}_{y-1}^{(x,m)} < n^{1/4} } > \frac{\sqrt{n}}{K} \right)
&\leq \sup_{m\geq 0} P\left( \sum_{i\leq K\sqrt{n}} \ind{\hat{U}_i < n^{1/4}} > \frac{\sqrt{n}}{2K} \, \biggl| \, \hat{U}_0 = m \right) \\
&\quad + \sup_{m\geq 1} P\left( \sum_{i\leq \s_0^U} \ind{U_i < n^{1/4}} > \frac{\sqrt{n}}{2K} \, \biggl| \, U_0 = m  \right),
\end{align*}
and Lemmas \ref{lem:BLPsmall1} and \ref{lem:BLPsmall2} imply that both probabilities on the right decay faster than any polynomial in $n$. 
It follows from this and Corollary \ref{range} and Lemma \ref{lt} that 
\begin{equation}\label{PGnc}
 \lim_{K\ra\infty} \limsup_{n\ra\infty} P(G_{n,K,t}^c) = 0. 
\end{equation}
If the event $G_{n,K,t}$ occurs, then for any $k\leq nt$ there is an $x \in [-K\sqrt{n},K\sqrt{n}]$ and $0\leq m\leq K\sqrt{n}$ such that $\l_{x,m}=k$. 
Moreover, on the event $G_{n,K,t} \cap \{ \l_{x,m} \leq nt \}$ for fixed $|x|,m\leq K\sqrt{n}$ we have 
\begin{align*}
 \left| \sum_{y=x+1}^{\mathcal{M}^{(x,m)}} (\rho_y^{(x,m)} - \rho) \right|
&\leq C \sum_{y=x+1}^{\mathcal{M}^{(x,m)}} \exp\left\{ -c \mathcal{E}_{y-1}^{(x,m)} \right\} \\
&\leq C (2K\sqrt{n}+1) e^{-c n^{1/4}} + C \sum_{y=x+1}^{\mathcal{M}^{(x,m)}} \ind{ \mathcal{E}_{y-1}^{(x,m)} < n^{1/4}} \\
&\leq C (2K\sqrt{n}+1) e^{-c n^{1/4}} + \frac{C \sqrt{n}}{K}. 
\end{align*}
If $K$ is chosen large enough so that $K>2C/\e$ then for $n$ sufficiently large the last line above is less than $\e\sqrt{n}/2$. 
Therefore, if $K > 2C/\e$ we have for all sufficiently large $n$ that
\begin{align}
 &P\left( \sup_{k\leq nt} |C_k^+ - \rho(M_k-X_k)| \geq \e \sqrt{n} \right) \nonumber \\
 &\leq P(G_{n,K,t}^c) + K^2 n \sup_{\substack{|x|\leq K\sqrt{n}\\ m\leq K\sqrt{n}}} P\left( \left| \sum_{y=x+1}^{\mathcal{M}^{(x,m)}} \left( \Delta_y^{(x,m)} - \rho \right) \right| \geq \e \sqrt{n}, \, G_{n,K,t}, \, \l_{x,m} \leq nt \right)  \nonumber \\
 &\leq P(G_{n,K,t}^c) + K^2 n \sup_{\substack{|x|\leq K\sqrt{n}\\ m\leq K\sqrt{n}}} P\left( \left| \sum_{y=x+1}^{\mathcal{M}^{(x,m)}} \left( \Delta_y^{(x,m)} - \rho_y^{(x,m)} \right) \right| \geq \frac{\e \sqrt{n}}{2} , \, G_{n,K,t}, \, \l_{x,m} \leq nt \right)  \nonumber \\
 &\leq P(G_{n,K,t}^c) + K^2 n \sup_{\substack{|x|\leq K\sqrt{n}\\ m\leq K\sqrt{n}}} P\left( \sup_{k\leq 2K\sqrt{n}} \left| \sum_{y=x+1}^{x+k} \left( \Delta_y^{(x,m)} - \rho_y^{(x,m)} \right) \right| \geq \frac{\e \sqrt{n}}{2} \right). \label{MGdiff} 
\end{align}
Since the random variables $\Delta_y^{(x,m)}$ only take values from the finite set  $\{ \sum_{j=1}^\ell (2p_j-1) , \, \ell=1,2,\ldots,N\}$,  the sum inside the probability in \eqref{MGdiff} is a martingale with bounded increments. 
Therefore, it follows from Azuma's inequality that the last probability in \eqref{MGdiff} is bounded above by $2 \exp\{-\frac{C \e^2 \sqrt{n}}{K} \}$ for some constant $C>0$ that does not depend on $x$ or $m$, and so the second term in \eqref{MGdiff} vanishes as $n\ra\infty$ for any fixed $K$. 
Finally, recalling \eqref{PGnc} the limit in \eqref{Cplus} follows by taking $n\ra\infty$ and then $K\ra\infty$ in \eqref{MGdiff}. 
\end{proof}

\textbf{Step 3: Tightness.}
The next step in the proof of Theorem \ref{thm:pbm} is to prove tightness for the 
random walk under diffusive scaling.

\begin{lem}\label{lem:Xtight}
The sequence of processes $X_{\fl{n\cdot}}/\sqrt{n}$, $n\ge 0$,
is tight in the space $D([0,\infty))$ of c\'adl\'ag paths equipped with the Skorokhod($J_1$) topology. Moreover, any subsequential limiting distribution is concentrated on continuous paths. 
\end{lem}

\begin{proof}
For the proof of this lemma, and also for step 4 of the proof, it will be helpful to use a slightly different decomposition of the walk instead of $X_n = B_n + C_n$. 
Lemma \ref{lem:Cnapprox} implies that $X_n$ can be approximated as
\[
 X_n \approx B_n + \rho M_n + \tilde\rho I_n - (\rho + \tilde\rho)X_n = B_n + \rho M_n + \tilde\rho I_n -  \left( \frac{\nu}{2}-1 \right) X_n, 
\]
where we used the identity \eqref{rhosums} in the last equality. 
The disadvantage of this decomposition is that $X_n$ appears both on the left and the right above. To account for this, we define $D_n = C_n + (\frac{\nu}{2}-1) X_n$ so that $X_n = B_n + D_n - (\frac{\nu}{2}-1) X_n$, or equivalently, 
\begin{equation}\label{XnDecomp}
 X_n = \frac{2}{\nu} B_n + \frac{2}{\nu} D_n.
\end{equation}
The representation \eqref{XnDecomp} is helpful since we know by Lemma \ref{lem:BnBM} that the first term on the right converges to Brownian motion and by Lemma \ref{lem:Cnapprox} that the second term on the right is approximated by a linear combination of the running maximum and running minimum of the walk. We will use these facts to reduce the proof of Lemma \ref{lem:Xtight} to proving tightness for the (rescaled) running maximum and running minimum processes.

The conclusions of Lemma \ref{lem:Xtight} will follow if we can show 
\begin{equation}\label{ERWtight}
 \lim_{\d \ra 0} \limsup_{n\ra\infty} P\left( \sup_{\substack{k,\ell \leq nt \\ |k-\ell|\leq n \d}} |X_k - X_\ell| \geq \e \sqrt{n} \right) = 0, \quad \forall \e>0, \, t<\infty. 
\end{equation}
Indeed, it follows from \eqref{ERWtight} and Corollary \ref{range} that the sequence $X_{\fl{n\cdot}}/\sqrt{n}$,  $n\ge 0$, is tight (see \cite[Theorem 16.8]{bCOPM}). 
Moreover, since the rescaled random walk has jumps of size $\pm 1/\sqrt{n} \xrightarrow{n\ra\infty} 0$ it follows that any subsequential limit of $X_{\fl{n\cdot}}/\sqrt{n}$ in the space $D([0,t])$ is concentrated on continuous paths (see \cite[Theorem 13.4]{bCOPM}). 
Using the decomposition \eqref{XnDecomp}, to prove \eqref{ERWtight} it will be enough to show 
\begin{equation}\label{Bntight}
 \lim_{\d \ra 0} \limsup_{n\ra\infty} P\left( \sup_{\substack{k,\ell \leq nt \\ |k-\ell|\leq n \d}} |B_k - B_\ell| \geq \e \sqrt{n} \right) = 0, \quad \forall \e>0, \, t < \infty, 
\end{equation}
and
\begin{equation}\label{Dntight}
 \lim_{\d \ra 0} \limsup_{n\ra\infty} P\left( \sup_{\substack{k,\ell \leq nt \\ |k-\ell|\leq n \d}} |D_k - D_\ell| \geq \e \sqrt{n} \right) = 0, \quad \forall \e>0, \, t < \infty. 
\end{equation}
The limit in \eqref{Bntight} follows from Lemma \ref{lem:BnBM}, and to prove \eqref{Dntight} note that 
\[
 \sup_{\substack{k,\ell \leq nt \\ |k-\ell|\leq n \d}} |D_k-D_\ell|
\leq 2 \sup_{k\leq nt} |D_k - \rho M_k - \tilde\rho I_k| + \sup_{\substack{k,\ell \leq nt \\ |k-\ell|\leq n \d}} |\rho| |M_k-M_\ell| + \sup_{\substack{ k,\ell \leq nt \\ |k-\ell| \leq n\d}} |\tilde\rho| |I_k-I_\ell|
\]
so that \eqref{Dntight} follows from Corollary \ref{cor:MnIntight} and Lemma \ref{lem:Cnapprox}. 
\end{proof}

\textbf{Step 4: Convergence to perturbed Brownian motion.}
Finally, we collect the results from the first three steps to prove that the rescaled path of the ERW converges in distribution to a $(\th,\tilde\th)$-perturbed Brownian motion. 
We begin by introducing the following notation for the rescaled process versions of $X_n$, $B_n$ and $D_n$, 
\[
\mathfrak{X}_n(t) = \frac{X_{\fl{nt}}}{\sqrt{ \frac{2}{\nu} n}}, 
\quad
\mathfrak{B}_n(t) = \frac{B_{\fl{nt}}}{\sqrt{\frac{\nu}{2} n}}, 
\quad\text{and}\quad
\mathfrak{D}_n(t) = \frac{D_{\fl{nt}}}{ \sqrt{\frac{\nu}{2} n}}, 
\]
so that \eqref{XnDecomp} implies that $\mathfrak{X}_n = \mathfrak{B}_n+\mathfrak{D}_n$. 
Note that Lemmas \ref{lem:BnBM} and \ref{lem:Xtight} imply that the joint sequence $(\mathfrak{X}_n, \mathfrak{B}_n, \mathfrak{D}_n)$ is a tight sequence in $D([0,\infty))^3$ such that any subsequence that converges in distribution is concentrated on $C([0,\infty))^3$. 
(Note that the argument in the proof of Lemma \ref{lem:Xtight} also shows that the sequence of paths $\mathfrak{D}_n$ is tight with any subsequential weak limits concentrated on continuous paths.)

Now, let $\Psi:D([0,\infty)) \ra D([0,\infty))$ be the mapping defined by 
\[
 \Psi(x)(t) = \th \, \sup_{s\leq t} x(s) + \tilde\th \, \inf_{s \leq t} x(s). 
\]
Note that since $\theta = (2\rho)/\nu$ and $\tilde\th = (2\tilde\rho)/\nu$ then Lemma \eqref{lem:Cnapprox} is equivalent to the statement that 
\begin{equation}\label{Dnapprox}
 \lim_{n\ra\infty} P\left( \sup_{s\leq t} \left| \mathfrak{D}_n(s) - \Psi(\mathfrak{X}_n)(s) \right| \geq \e \right) = 0, \quad \forall \e>0, \, t<\infty. 
\end{equation}
Since the function $\Psi$ is continuous on the subset $C([0,\infty))$ of continuous functions, it follows from \eqref{Dnapprox}, the continuous mapping theorem, and Lemmas \ref{lem:BnBM} and \ref{lem:Xtight} that if $n_k$ is a subsequence on which the joint sequence $(\mathfrak{X}_n, \mathfrak{B}_n, \mathfrak{D}_n)$ converges in distribution it must converge to a joint process of the form $(\mathfrak{X}, \mathfrak{B}, \Psi(\mathfrak{X}))$, where $\mathfrak{X}$ is a continuous process and $\mathfrak{B}$ is a standard Brownian motion. However, since $\mathfrak{X}_n= \mathfrak{B}_n + \mathfrak{D}_n$ the limit process must also satisfy $\mathfrak{X} = \mathfrak{B} + \Psi(\mathfrak{X})$; that is, 
\[
 \mathfrak{X}(t) = \mathfrak{B}(t) + \th \sup_{s\leq t} \mathfrak{X}(s) + \tilde\th \inf_{s\leq t} \mathfrak{X}(s), \qquad \forall t\geq 0, 
\]
and so $\mathfrak{X}$ must be a $(\th,\tilde\th)$-perturbed Brownian motion. 
Note that the above argument shows that any subsequential limit $\mathfrak{X}_n$ is a $(\th,\tilde\th)$-perturbed Brownian motion, and since the sequence $\mathfrak{X}_n$ is tight it follows that $\mathfrak{X}_n$ converges in distribution to a $(\th,\tilde\th)$-perturbed Brownian motion. 

\section{Markovian cookie stacks}\label{sec:Markov}


The main results for this paper (Theorems \ref{thm:boundary} and
\ref{thm:pbm}) are stated for recurrent ERW with periodic cookie
stacks. However, many of the results in the present paper also hold
for recurrent ERW in the more general model of Markovian cookie stacks
which was introduced and studied in \cite{kpERWMC}.  For
Markovian cookie stack environments, the cookie environment
$\w = \{\w_x\}_{x\in\Z}$ is spatially i.i.d., but the cookie sequence
at each site $\w_x = \{\w_x(j)\}_{j\geq 1}$ comes from the realization
of a finite state Markov chain. For this model the
  analogs of parameters $\theta,\ \tilde{\theta},\ \nu$ have been
  explicitly computed in \cite{kpERWMC}.\footnote{Parameters
    $\delta,\ \tilde{\delta}$ in \cite{kpERWMC} correspond to
    $\theta,\ \tilde{\theta}$ of the current paper.} Nearly all of
the results in the present paper can be adapted to the case of
recurrent ERW in Markovian cookie stacks with little or no changes in
the proofs. In particular, since the proof of the functional limit
laws in the boundary cases (Theorem \ref{thm:boundary}) depends only
on the tail asymptotics for the BLP in Theorem \ref{thm:BLPtail}, and
since these tail asymptotics were proved for the model of Markovian
cookie stacks in \cite{kpERWMC}, then Theorem \ref{thm:boundary} also
holds for ERW in Markovian cookie stacks.

\begin{figure}[ht]
 \includegraphics[scale=0.4]{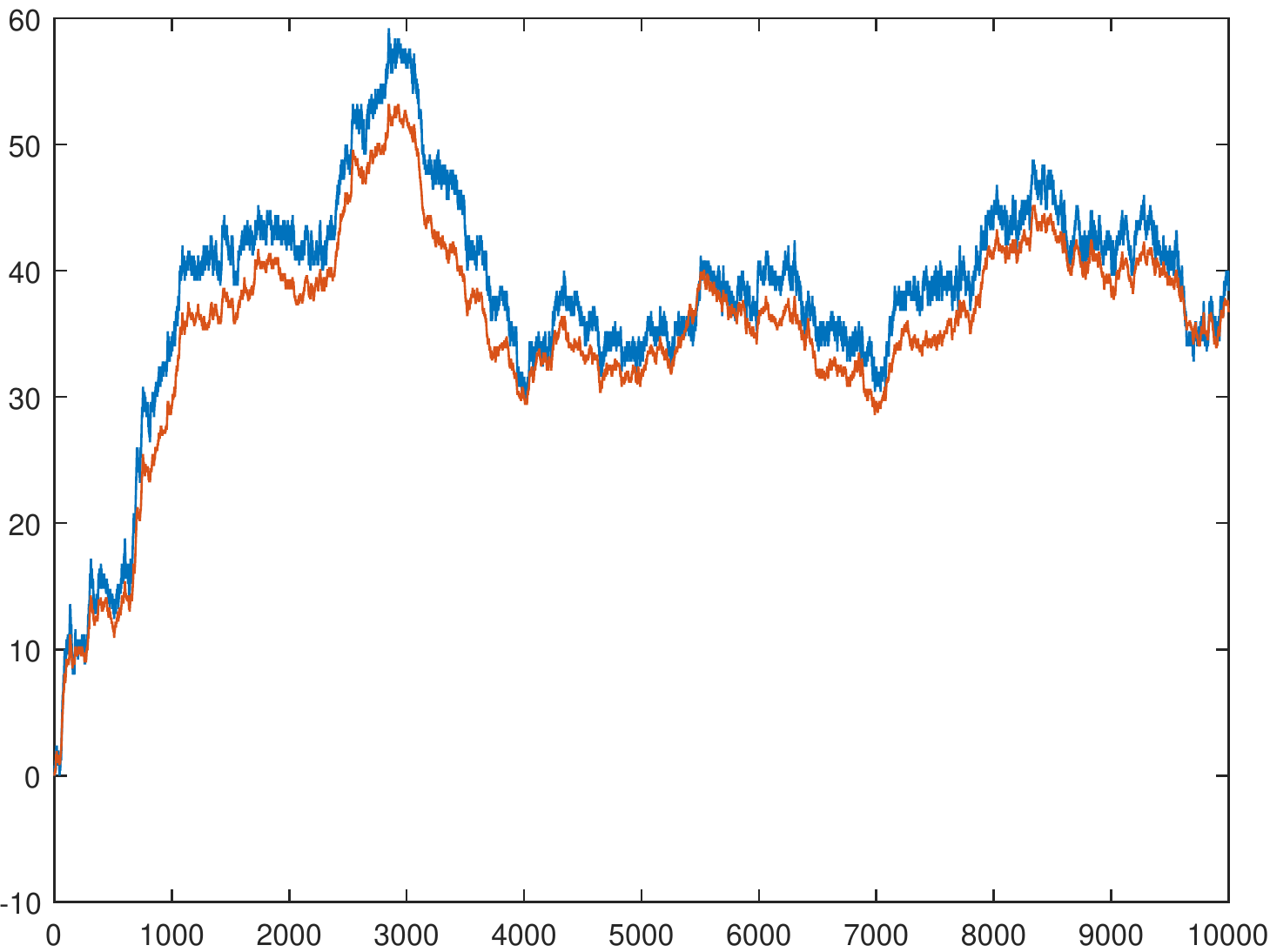}
 \includegraphics[scale=0.4]{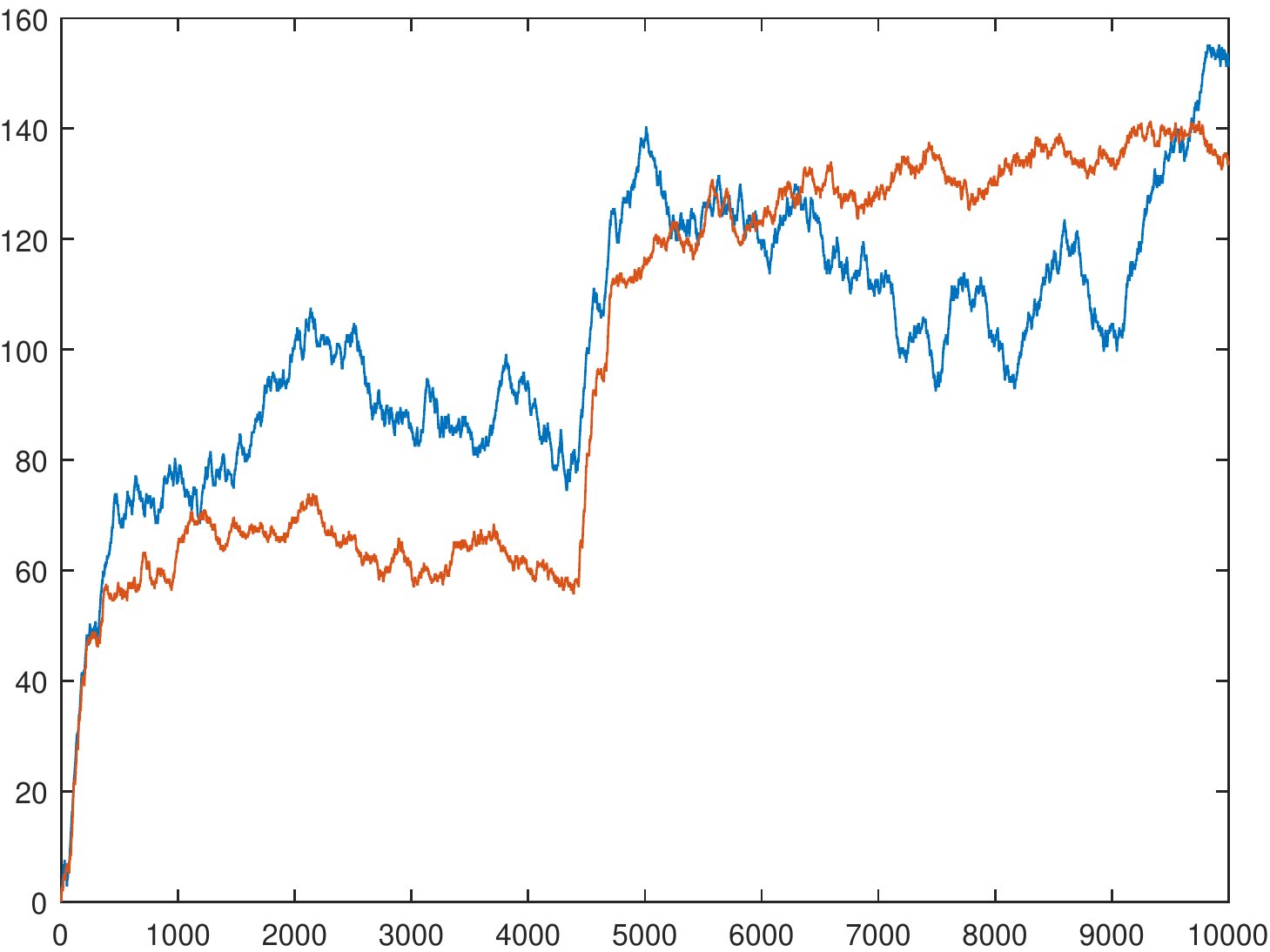}
\caption[Simulations of comparisons of $C_n$ with $\rho(M_n-X_n)+\tilde\rho(I_n-X_n)$]{The above plots both give the plot of the process $C_n$ (in blue) compared with the plot of $\rho(M_n-X_n) + \tilde\rho(I_n-X_n)$ (in red). The plot on the left comes from a simulation of an ERW with periodic cookie stacks of the form $(0.7,0.3,0.7,0.3,\ldots)$ at each site. The plot on the right comes from a simulation of an ERW with Markovian cookie stacks where the sequence $\{\w_x(j)\}_{j\geq 1}$ at each site is a Markov chain taking values in $\{0.7,0.3\}$, with transition matrix 
$K = \left(\begin{smallmatrix} 0.75 & 0.25 \\ 0.25 & 0.75 \end{smallmatrix}\right)$ 
and initial value $\w_x(1) = 0.7$.  }\label{fig:Cnplots}
\end{figure}

The only part of the current paper that doesn't generalize to
Markovian cookie stacks is Lemma \ref{lem:Cnapprox} which controls the
error in the approximation of $C_n$ by
$\rho(M_n-X_n) + \tilde\rho(I_n-X_n)$. In fact, in the general case of
Markovian cookie stacks a heuristic argument suggests that the difference
$C_n - \rho(M_n-X_n) - \tilde\rho(I_n-X_n)$ is of the order
$\sqrt{n}$, whereas Lemma \ref{lem:Cnapprox} shows that it is
$o(\sqrt{n})$ for periodic cookie stacks. 
(See Figure \ref{fig:Cnplots} for simulations which support this claim.)
We remark, that the only place where the proof of Lemma \ref{lem:Cnapprox} breaks down for Markovian cookie stacks is in the application of Azuma's inequality to bound the second probability in \eqref{MGdiff}. 
Thus, Lemma \ref{lem:Cnapprox} and therefore also Theorem \ref{thm:pbm} hold for more general recurrent ERW with the property that the partial sums $\sum_{j=1}^n (2\w_x(j) - 1)$ are uniformly bounded in $x$ and $n$. 

In spite of the fact that Lemma \ref{lem:Cnapprox} does not hold for the more general Markovian cookie stack model, other results suggest that it may be the case that the scaling limit for recurrent ERW in the non-boundary cases are still perturbed Brownian motions. 
\begin{conj}\label{con:MCS}
  If $X_n$ is an ERW in a cookie environment with Markovian cookie
  stacks and $\max\{\th,\tilde\th\} < 1$ , then for some $a>0$ the
  sequence of processes
  $\left\{\frac{X_{\fl{nt}}}{a\sqrt{n}} \right\}_{t\ge 0}$ converges
  in distribution to a $(\th,\tilde\th)$-perturbed Brownian motion.
\end{conj}
Some evidence for Conjecture \ref{con:MCS} is provided by the diffusion approximations for the BLPs in Theorem \ref{thm:da} which were already proved in the case of Markovian cookie stacks in \cite{kpERWMC}. Due to the connection of the BLPs with directed edge local times of the ERW, these diffusion approximations are consistent with Ray-Knight theorems for the local times of perturbed Brownian motion that were proved in \cite{cpyBetaPBM}. 

\begin{rem}
In the process of finishing the current paper, we learned of a recent paper \cite{HLSH16} which proves the convergence to perturbed Brownian motion for a special case of ERW in Markovian cookie stacks. This paper  considers the case where the cookie sequence at each site $\{\omega_x(j)\}_{j\geq 1}$ is a Markov chain on $\{0,1\}$ with transition probability $\begin{pmatrix} p & 1-p \\ 1-p & p \end{pmatrix}$. 
Note that for this model the path of the random walk is deterministic once the cookie environment $\w$ is fixed since $\w_x(j) \in \{0,1\}$ for all $x\in \Z$, $j\geq 1$.  
Convergence to perturbed-Brownian motion remains an open problem for the more general model considered in Conjecture \ref{con:MCS}.
\end{rem}




\bibliographystyle{alpha}
\bibliography{CookieRW}



\end{document}